\documentclass[11pt]{article}
\usepackage{amsfonts}
\usepackage[mathscr]{eucal}
\usepackage{mathrsfs}
\usepackage{amsmath}
\usepackage{amssymb}
\usepackage{amsmath}
\usepackage{amsthm}
\usepackage{mathrsfs}
\usepackage{graphicx}

\usepackage{url,hyperref,lineno,microtype,subcaption}
\usepackage[onehalfspacing]{setspace}

\newtheorem{theorem}{Theorem}[section]
\newtheorem{corollary}{Corollary}
\newtheorem{lemma}[theorem]{Lemma}

\newtheorem{definition}[theorem]{Definition}
\newtheorem{remark}{Remark}

\newcommand{\Ga}{{\Gamma}}

\newcommand{\Om}{{\Omega}}

\newcommand{\be}{{\beta}}
\newcommand{\ga}{{\gamma}}
\newcommand{\de}{{\delta}}
\newcommand{\ep}{\varepsilon}

\newcommand{\la}{\lambda}
\newcommand{\si}{{\sigma}}

\newcommand{\vph}{{\varphi}}
\newcommand{\om}{{\omega}}

\newcommand{\F}{{\mathbb F}}

\newcommand{\R}{{\mathbb R}}

\newcommand{\cB}{\mathcal{B}}

\newcommand{\cF}{\mathcal{F}}

\newcommand{\cE}{\mathcal{E}}

\newcommand{\cP}{\mathcal{P}}

\newcommand{\cT}{\mathcal{T}}

\newcommand{\cK}{\mathcal{K}}

\newcommand{\cW}{\mathcal{W}}

\newcommand{\g}{{\nabla}}
\newcommand{\pd}{\partial}
\newcommand{\intl}{\int\limits}
\newcommand{\arr}{\rightarrow}

\newcommand{\oPhi}{\overline{\Phi}}

\newcommand{\oom}{\overline{\om}}
\newcommand{\ovph}{\overline{\varphi}}
\newcommand{\ou}{\overline{u}}
\newcommand{\ov}{\overline{v}}
\newcommand{\ow}{\overline{w}}

	\title {Qualitative properties of solutions to  a  nonlinear transmission problem for an elastic Bresse beam}

\author{Tamara Fastovska\,$^{1,2,*}$, Dirk Langemann\,$^{3}$ and Iryna Ryzhkova\,$^{1}$} 

\begin{document}
		\maketitle
		{$^{1}$Department of Mathematics and Computer Science, V.N. Karazin Kharkiv National University, Kharkiv, Ukraine \\
			$^{2}$ Institut für Mathematik, Humboldt-Universität zu Berlin, Berlin, Germany\\
			$^{3}$Institut für Partielle Differentialgleichungen,Technische Universität Braunschweig, Braunschweig, Germany}\\
	

	\begin{abstract}
		\section{}
		We consider a nonlinear transmission problem for a Bresse beam, which consists of two parts, damped and undamped. The mechanical damping in the damped part is present  in the shear angle equation only, and the damped part may be of arbitrary positive length. We prove well-posedness of the corresponding PDE system in  energy space and establish existence of a regular global attractor under certain conditions on nonlinearities and coefficients of the damped part only. Moreover, we study singular limits of the problem  when $l\arr 0$ or $l\arr 0$ simultaneously with $k_i\arr +\infty$ and perform numerical modelling for these processes.
		
		\tiny
		 \section{Keywords:} Bresse beam, transmission problem, global attractor, singular limit
	\end{abstract}

	\section{Introduction}
\begin{figure}[h]
		\centering
		\includegraphics[width=0.45\textwidth]{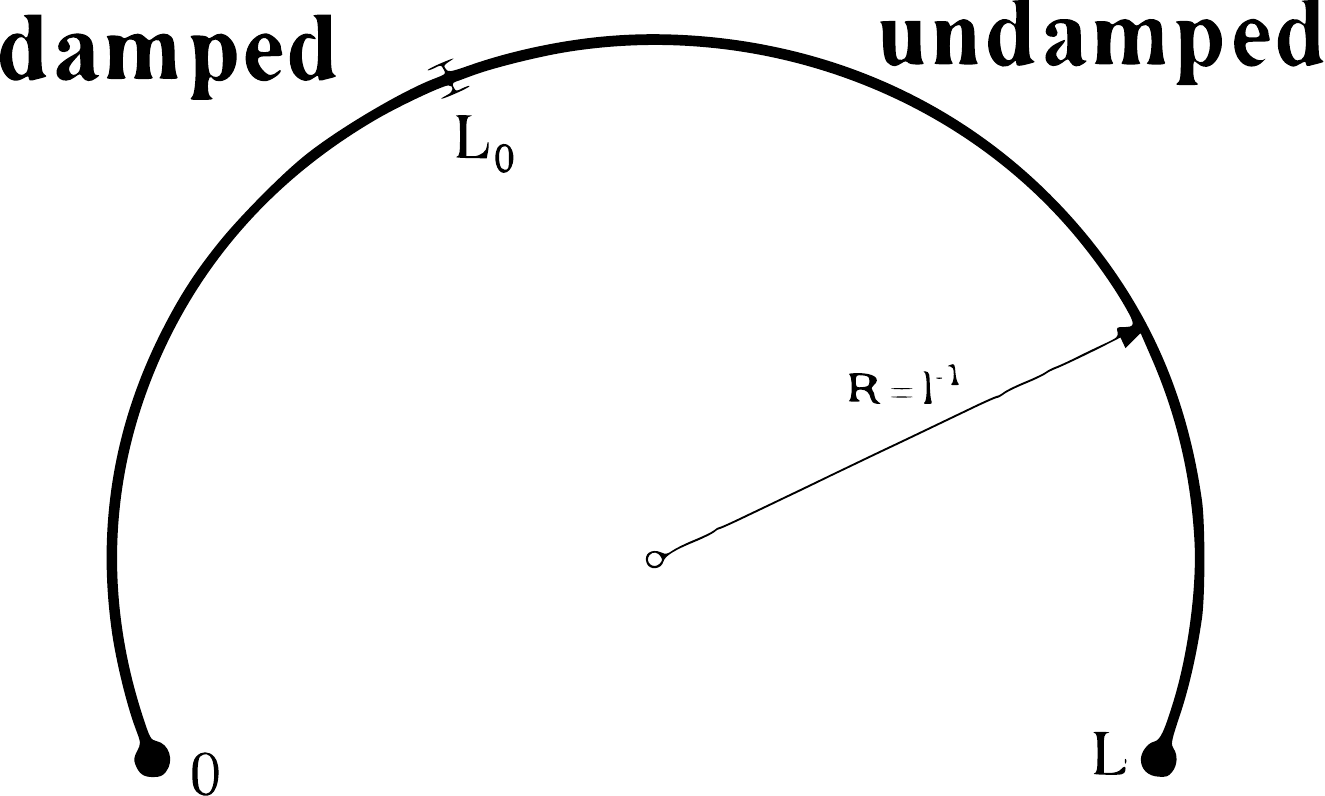}
		\caption{Composite Bresse beam.} \label{FigBeam}
	\end{figure}

	In this paper we consider a contact problem for the Bresse beam. Originally the mathematical model for homogeneous Bresse beams was derived in \cite{Bre1859}. We use the variant of the model described in \cite[Ch. 3]{LagLeu1994}. Let the whole beam occupies a part of a circle of length $L$ and has the curvature $l=R^{-1}$. We consider the beam as a one-dimensional object and measure the coordinate $x$ along the beam. Thus, we say that the coordinate $x$ changes within the interval $(0,L)$.  The parts of the beam occupying the intervals $(0,L_0)$ and  $(L_0,L)$ consist of different materials. The part lying in the interval $(0,L_0)$ is partially subjected to a structural damping (see Figure \ref{FigBeam}).
	The Bresse system describes evolution of three quantities: transversal displacement, longitudinal displacement and shear angle variation. We denote by $\vph$, $\psi$, and $\om$  the transversal displacement, the shear angle variation, and  the longitudinal displacement of the left part of the beam lying in $(0,L_0)$. Analogously, we denote by  $u$, $v$, and $w$  the transversal displacement,  the shear angle variation,  and the longitudinal displacement of the right part of the beam occupying the interval $(L_0,L)$.  We assume the presence of the mechanical dissipation in the equation for the shear angle variation for the left part of the beam. We also assume that both ends of the beam are clamped.
	Nonlinear oscillations of the composite beam can be described by the following  system of equations
	\begin{align}
		& \rho_1\vph_{tt}-k_1(\vph_x+\psi+l\om)_x - l\si_1(\om_x-l\vph) +f_1(\vph, \psi, \om)=p_1(x,t), \label{Eq1}\\
		& \be_1\psi_{tt} -\la_1 \psi_{xx} +k_1(\vph_x+\psi+l\om) +\ga(\psi_t) +h_1(\vph, \psi, \om)=r_1(x,t),\;x\in (0,L_0), t>0,\label{Eq2}\\
		& \rho_1\om_{tt}- \si_1(\om_x-l\vph) _x+lk_1(\vph_x+\psi+l\om)+g_1(\vph, \psi, \om)=q_1(x,t), \label{Eq3}
	\end{align}
	and
	\begin{align}
		& \rho_2u_{tt}-k_2(u_x+v+lw)_x - l\si_2(w_x-lu) +f_2(u, v, w)=p_2(x,t), \label{Eq4}\\
		& \be_2v_{tt} -\la_2 v_{xx} +k_2(u_x+v+lw) +h_2(u, v, w)=r_2(x,t),\qquad \qquad x\in (L_0,L), t>0,\label{Eq5}\\
		& \rho_2w_{tt}- \si_2(w_x-lu)_x+lk_2(u_x+v+lw)+g_2(u, v, w)=q_2(x,t), \label{Eq6}
	\end{align}
	where $\rho_j,\;\be_j, \; k_j, \; \si_j,\; \la_j$ are positive parameters,   $f_j, \; g_j, \; h_j:\R^3\arr\R$ are nonlinear feedbacks,  $p_j, \; q_j, \; r_j:(0,L)\times\R^3\arr\R$ are known external loads, $\gamma:\R\arr\R$ is a nonlinear damping.
	The system is subjected to the Dirichlet boundary conditions
	\begin{equation} \label{BC}
		\vph(0,t)=u(L,t)=0,  \quad \psi(0,t)=v(L,t)=0, \quad \om(0,t)=w(L,t)=0,
	\end{equation}
	the  transmission conditions
	\begin{align}
		& \vph(L_0,t)=u(L_0,t), \quad  \psi(L_0,t)=v(L_0,t),  \quad  \om(L_0,t)=w(L_0,t), \label{TC1} \\
		& k_1(\vph_x+\psi+l\om)(L_0,t)=k_2(u_x+v+lw)(L_0,t), \\
		& \la_1 \psi_{x}(L_0,t)= \la_2 v_{x}(L_0,t),\\
		& \si_1(\om_x-l\vph)(L_0,t)=\si_2(w_x-lu)(L_0,t), \label{TC4}
	\end{align}
	and supplemented with the initial conditions
	\begin{align}
		&\vph(x,0)=\vph_0(x),\quad \psi(x,0)=\psi_0(x),\quad \om(x,0)=\om_0(x),\\
		&\vph_t(x,0)=\vph_1(x),\quad \psi_t(x,0)=\psi_1(x),\quad \om_t(x,0)=\om_1(x),
	\end{align}
	\begin{align}
		&u(x,0)=u_0(x),\quad v(x,0)=v_0(x),\quad w(x,0)=w_0(x),\\
		&u_t(x,0)=u_1(x),\quad v_t(x,0)=v_1(x),\quad w_t(x,0)=w_1(x).\label{IC}
	\end{align}
	One can observe patterns in the problem which appear to have physical meaning:
	\begin{align*}
		& Q_i(\xi, \zeta, \eta)=k_i(\xi_x+\zeta +l\eta) \mbox{  are  shear forces},\\
		& N_i(\xi, \zeta, \eta)=\si_i(\eta_x-l\xi) \mbox{ are axial forces},\\
		& M_i(\xi, \zeta, \eta)=\la_i\zeta_x \mbox{ are bending moments}
	\end{align*}
	for damped ($i=1)$ and undamped ($i=2$) parts respectively.
	Later we will use them to rewrite the problem in a compact and physically natural form.
	
	The paper is devoted to the well-posedness and long-time behaviour of the system \eqref{Eq1}-\eqref{IC}. Our main goal is to establish conditions under which the assumed amount of dissipation is sufficient to guarantee the existence of a global attractor.
	
	The paper is organized as follows. In Section 2 we represent functional spaces and pose the problem in an abstract form. In Section 3 we prove that the problem is well-posed and possesses strong solutions provided nonlinearities and initial data are smooth enough.  Section 4 is devoted  to the main result on the existence of a compact attractor. The  nature of dissipation prevents us from proving dissipativity explicitly, thus we show that the corresponding dynamical system is of gradient structure and asymptotically smooth. We establish the unique continuation property by means of the observability inequality obtained in \cite{TriYao2002} to prove the gradient property.  The compensated compactness approach is used to prove the asymptotic smoothness. In Section 5 we show that solutions to \eqref{Eq1}-\eqref{IC} tend to solutions to a transmission problem for the Timoshenko beam when $l\arr 0$ and to solutions to a transmission problem for the Euler-Bernoulli beam when $l\arr 0$ and $k_i \arr\infty$ as well as perform numerical modelling of these singular limits.

	\section{Preliminaries and Abstract formulation}
	\subsection{Spaces and notations}
	Let us denote
	\begin{equation*}
		\Phi^1=(\vph, \psi, \om), \quad\Phi^2=(u,v,w), \quad \Phi=(\Phi^1,\Phi^2).
	\end{equation*}
	Thus, $\Phi$ is a six-dimensional vector of functions. Analogously,
	\begin{align*}
		& F_j=(f_j,g_j, h_j): \R^3\arr\R^3, \quad F=(F_1,F_2), \\
		& P_j=(p_j,q_j, r_j): [(0,L)\times \R_+]^3\arr\R^3, \quad P=(P_1,P_2), \\
		& R_j=diag\{\rho_j,\be_j,\rho_j\}, \quad R=diag\{\rho_1,\be_1,\rho_1,\rho_2,\be_2,\rho_2\}, \\
		& \Ga(\Phi_t)=(0,\ga(\psi_{t}),0,0,0,0),
	\end{align*}
	where $j=1,2$. The static linear part of the equation system  can be formally rewritten as
	\begin{equation}\label{oper}
		A\Phi=\left(
		\begin{aligned}
			& -\pd_x Q_1(\Phi^1)-lN_1(\Phi^1)\\
			& -\pd_x M_1(\Phi^1)+Q_1(\Phi^1)\\
			& -\pd_x N_1(\Phi^1)+lQ_1(\Phi^1) \\
			& -\pd_x Q_2(\Phi^2)-lN_2(\Phi^2)\\
			& -\pd_x M_2(\Phi^2)+Q_2(\Phi^2) \\
			& -\pd_x N_2(\Phi^2)+lQ_2(\Phi^2)
		\end{aligned}
		\right).
	\end{equation}
	Then transmission conditions \eqref{TC1}-\eqref{TC4} can be written as follows
	\begin{align*}
		&\Phi^1(L_0,t)=\Phi^2(L_0,t), \\
		& Q_1(\Phi^1(L_0,t))=Q_2(\Phi^2(L_0,t)),\\
		& M_1(\Phi^1(L_0,t))=M_2(\Phi^2(L_0,t)),\\
		& N_1(\Phi^1(L_0,t))=N_2(\Phi^2(L_0,t)).
	\end{align*}
	Throughout  the paper we use the notation $||\cdot||$ for the $L^2$-norm of a function and $(\cdot,\cdot)$ for the $L^2$-inner product. In these notations we skip  the domain, on which functions are defined. We adopt the notation $||\cdot||_{L^2(\Omega)}$ only when domain is not evident.  We also use the same notations  $||\cdot||$ and  $(\cdot,\cdot)$ for $[L^2(\Om)]^3$.\\
	To write our problem in an abstract form we introduce the following spaces. For the velocities of the displacements we use the space
	\begin{equation*}
		H_v=\{\Phi=(\Phi^1,\Phi^2):\; \Phi^1\in [L^2(0,L_0)]^3, \;\Phi^2\in [L^2(L_0,L)]^3 \}
	\end{equation*}
	with the norm
	\begin{equation*}
		||\Phi||^2_{H_v}=||\Phi||^2_v=\sum_{j=1}^{2}||\sqrt{R_j}\Phi^j||^2,
	\end{equation*}
	which is equivalent to the standard $L^2$-norm.\\
	For the beam displacements we use the space
	\begin{align*}
		H_d=\left\{\Phi\in H_v:\; \right.&\Phi^1\in [H^1(0,L_0)]^3, \;\Phi^2\in [H^1(L_0,L)]^3,  \\
		&\left.\Phi^1(0,t)=\Phi^2(L,t)=0,\; \Phi^1(L_0,t)=\Phi^2(L_0,t) \right\}
	\end{align*}
	with the norm
	\begin{equation*}
		||\Phi||^2_{H_d}=||\Phi||^2_d=\sum_{j=1}^{2}\left(||Q_j(\Phi^j)||^2+||N_j(\Phi^j)||^2+||M_j(\Phi^j)||^2 \right).
	\end{equation*}
	This norm is equivalent to the standard $H^1$-norm. Moreover, the equivalence constants can be chosen independent of $l$ for $l$ small enough (see \cite{MaMo2017}, Remark 2.1). If we  set
	\begin{equation*}
		\Psi(x)=\left\{
		\begin{aligned}
			&\Phi^1(x), \quad &x\in (0,L_0),\\
			&\Phi^2(x), \quad &x\in [L_0,L)
		\end{aligned}
		\right.
	\end{equation*}
	we see that there is isomorphism between $H_d$ and $[H^1_0(0,L)]^3$.
	\subsection{Abstract formulation}
	The operator $A:D(A)\subset H_v\arr H_v$ is defined by formula \eqref{oper}, where
	\begin{multline*}
		D(A)=\left\{\Phi\in H_d:\right. \Phi^1\in H^2(0,L_0), \; \Phi^2\in H^2(L_0,L), \;
		Q_1(\Phi^1(L_0,t))=Q_2(\Phi^2(L_0,t)), \\
		N_1(\Phi^1(L_0,t))=N_2(\Phi^2(L_0,t)),
		\left.M_1(\Phi^1(L_0,t))=M_2(\Phi^2(L_0,t))\right.\}
	\end{multline*}
	Arguing analogously to  Lemmas 1.1-1.3 from \cite{LiuWil2000} one can prove the following lemma.
	\begin{lemma}\label{lem:ASelfAdjont}
		The operator $A$ is positive  and self-adjoint. Moreover,
		\begin{equation} \label{AForm}
			\begin{aligned}
				(A^{1/2}\Phi, A^{1/2}B)=
				& \frac{1}{k_1} (Q_1(\Phi^1),Q_1(B^1)) + \frac{1}{\si_1} (N_1(\Phi^1),N_1(B^1)) + \frac{1}{\la_1} (M_1(\Phi^1),M_1(B^1)) + \\
				& \frac{1}{k_2} (Q_2(\Phi^2),Q_2(B^2)) + \frac{1}{\si_2} (N_2(\Phi^2),N_2(B^2)) + \frac{1}{\la_2} (M_2(\Phi^2),M_2(B^2))
			\end{aligned}
		\end{equation}
		and $D(A^{1/2})=H_d\subset H_v$.
	\end{lemma}
	Thus,  we can rewrite equations \eqref{Eq1}-\eqref{Eq6} in the form
	\begin{equation} \label{AEq}
		R\Phi_{tt}+A\Phi+\Ga(\Phi_t) + F(\Phi)=P(x,t),
	\end{equation}
	boundary conditions \eqref{BC}  in the form
	\begin{equation}
		\Phi^1(0,t)=\Phi^2(L,t)=0, \label{ABC}
	\end{equation}
	and transmission conditions \eqref{TC1}-\eqref{TC4} can be written as
	\begin{align}
		&\Phi^1(L_0,t)=\Phi^2(L_0,t), \label{ATC1} \\
		& Q_1(\Phi^1(L_0,t))=Q_2(\Phi^2(L_0,t)),\\
		& M_1(\Phi^1(L_0,t))=M_2(\Phi^2(L_0,t)),\\
		& N_1(\Phi^1(L_0,t))=N_2(\Phi^2(L_0,t)). \label{ATC4}
	\end{align}
	Initial conditions have the form
	\begin{equation}
		\Phi(x,0)=\Phi_0(x), \qquad \Phi_t(x,0)=\Phi_1(x). \label{AIC}
	\end{equation}
	We use  $H=H_d\times H_v$ as a phase space.
	\section{Well-posedness}
	In this section we study strong, generalized and variational (weak) solutions to \eqref{AEq}-\eqref{AIC}.
	\begin{definition}
		$\Phi\in C(0,T;H_d)\bigcap C^1(0,T;H_v)$ such that $\Phi(x,0)=\Phi_0(x)$, $\Phi_t(x,0)=\Phi_1(x)$ is said to be  a strong solution to \eqref{AEq}-\eqref{AIC} if
		\begin{itemize}
			\item $\Phi(t)$ lies in $D(A)$ for almost all $t$;
			\item  $\Phi(t)$ is an absolutely continuous function with values in $H_d$  and $\Phi_t \in L_1(a,b;H_d)$ for $0<a<b<T$;
			\item $\Phi_t(t)$ is an absolutely continuous function with values in $H_v$  and $\Phi_{tt} \in L_1(a,b;H_v)$ for $0<a<b<T$;
			\item equation \eqref{AEq} is satisfied for almost all $t$.
		\end{itemize}
	\end{definition}
	\begin{definition}
		$\Phi\in C(0,T;H_d)\bigcap C^1(0,T;H_v)$ such that $\Phi(x,0)=\Phi_0(x)$, $\Phi_t(x,0)=\Phi_1(x)$ is said to be  a generalized solution to \eqref{AEq}-\eqref{AIC} if there exists a sequence of strong solutions $\Phi^{(n)}$ to \eqref{AEq}-\eqref{AIC} with the initial data $(\Phi_0^{(n)},\Phi_1^{(n)})$ and right hand side $P^{(n)}(x,t)$ such that
		\begin{equation*}
			\lim_{n\arr\infty} \max_{t\in[0,T]} \left(||\Phi^{(n)}(\cdot,t)-\Phi(\cdot,t)||_d  + ||\Phi_t^{(n)}(\cdot,t)-\Phi_t(\cdot,t)||_v \right) =0.
		\end{equation*}
	\end{definition}
	We also need a definition of a variational solution.
	We use six-dimensional vector-functions  $B=(B^1, B^2)$, $B^j=(\beta^j,\gamma^j,\delta^j)$ from the space
	\begin{equation*}
		F_T=\{B\in L^2(0,T;H_d), \; B_t\in L^2(0,T;H_v), B(T)=0\}
	\end{equation*}
	as test functions.
	\begin{definition}
		$\Phi$ is said to be a variational (weak) solution to \eqref{AEq}-\eqref{AIC} if
		\begin{itemize}
			\item $\Phi\in L^\infty(0,T;H_d), \; \Phi_t\in L^\infty(0,T;H_v)$;
			\item $\Phi$ satisfies the following variational equality for all $B\in F_T$
			\begin{multline}\label{VEq}
				-\intl_0^T  (R\Phi_t,B_t)(t)dt- (R\Phi_1, B(0)) + \int_0^T  (A^{1/2}\Phi, A^{1/2}B)(t)dt +\\
				\int_0^T  (\Gamma(\Phi_t), B)(t)dt +  \int_0^T  (F(\Phi), B)(t)dt - \int_0^T (P, B)(t)dt=0;
			\end{multline}
			\item $\Phi(x,0)=\Phi_0(x)$.
		\end{itemize}
	\end{definition}
	Now we state a well-posedness result for problem \eqref{AEq}-\eqref{AIC}.
	\begin{theorem}[Well-posedness] \label{th:WeakWP}
		Let
		\begin{align*}
			& f_i, \; g_i, \; h_i: \R^3 \arr \R \mbox{ are locally Lipschitz i.e.}  \\
			& |f_i(a)-f_i(b)|\le L(K)|a-b|, \quad \mbox{provided } |a|,|b|\le K; \tag{N1}\label{NLip}
		\end{align*}
		\begin{align*}
			& \mbox{there exist } \cF_i: \R^3\arr \R \mbox{ such that }  (f_i,h_i,g_i)=\g \cF_i; \\
			&\mbox{ there exists }  \delta>0   \mbox{ such that } \cF_j(a) \ge -\delta \mbox{ for all } a\in \R^3; \label{NBoundBelow} \tag{N2} \\
		\end{align*}
		\begin{equation*}
			P\in L^2(0,T;H_v);  \label{RSmooth}\tag{R1}
		\end{equation*}
		and the  nonlinear dissipation satisfies
		\begin{equation*}\label{DCont}
			\gamma \in C(\R) \mbox{ and non-decreasing }, \quad \gamma(0) = 0.  \tag{D1}
		\end{equation*}
		Then for every initial data $\Phi_0\in H_d, \Phi_1\in H_v$ and time interval $[0,T]$ there exists a unique generalized  solution to \eqref{AEq}-\eqref{AIC} with the following properties:
		\begin{itemize}
			\item every generalized solution is  variational;
			\item energy inequality 
			\begin{equation}\label{EE}
				\cE(T)+\int_0^T  (\ga(\psi_t),\psi_t)dt \le \cE(0)+ \int_0^T  (P(t),\Phi_t(t))dt
			\end{equation}
			holds, where
			\begin{equation*}
				\cE(t)=\frac 12 \left[||R^{1/2}\Phi_t(t)||^2+||A^{1/2}\Phi(t)||^2\right] + \intl_0^L  \cF(\Phi(x,t))dx
			\end{equation*}
			and
			\begin{equation*}
				\cF(\Phi(x,t))=\left\{
				\begin{aligned}
					& \cF_1(\vph(x,t),\psi(x,t),\om(x,t)), \quad & x\in (0,L_0),\\
					& \cF_2(u(x,t),v(x,t),w(x,t)), \quad & x\in (L_0,L).
				\end{aligned}
				\right.
			\end{equation*}
			\item If, additionally, $\Phi_0\in D(A)$, $\Phi_1\in H_d$ and
			\begin{equation*}
				\pd_t P(x,t)\in L_2(0,T;H_v)   \label{RAddSmooth}	\tag{R2}
			\end{equation*}
			then the generalized solution is also strong and satisfies the energy equality.
		\end{itemize}
	\end{theorem}
	\begin{proof}
		The proof essentially uses monotone operators theory. It  is rather standard by now (see, e.g., \cite{ChuEllLa2002}), so in some parts we give only references to corresponding arguments.  However, we give some details which demonstrate the peculiarity  of  1D problems.\\
		{\it Step 1. Abstract formulation.} We need to reformulate  problem \eqref{AEq}-\eqref{AIC} as a first order problem. Let us denote
		\begin{equation*}
			U=(\Phi,\Phi_t), \quad U_0=(\Phi_0, \Phi_1) \in H=H_d\times H_v,
		\end{equation*}
		\begin{equation*}
			\cT U=
			\begin{pmatrix}
				I & 0 \\
				0 & R^{-1}
			\end{pmatrix}
			\begin{pmatrix}
				0 & -I \\
				A & 0
			\end{pmatrix}
			U +
			\begin{pmatrix}
				0\\
				\Ga(\Phi_t)
			\end{pmatrix}.
		\end{equation*}
		Consequently, $D(\cT)=D(A)\times H_d\subset H$. In what follows  we use the notations
		\begin{equation*}
			\cB(U)=
			\begin{pmatrix}
				I & 0 \\
				0 & R^{-1}
			\end{pmatrix}
			\begin{pmatrix}
				0  \\
				F(\Phi)
			\end{pmatrix}, \quad
			\cP(x,t)=\begin{pmatrix}
				0  \\
				P(x,t)
			\end{pmatrix}.
		\end{equation*}
		Thus, we can rewrite problem \eqref{AEq}-\eqref{AIC} in the form
		\begin{equation}\label{FirstOrderForm}
			U_t+\cT U +\cB(U) = \cP, \quad U(0)=U_0\in H.
		\end{equation}
		{\it Step 2. Existence and uniqueness of a local solution.} Here we use Theorem~7.2 from \cite{ChuEllLa2002}.  For the reader's convenience we formulate it below.
		\begin{theorem}[\cite{ChuEllLa2002}]\label{th:AbstrEx}
			Consider the initial value problem
			\begin{equation}\label{AbstrIVP}
				U_t+\cT U +B(U) = f, \quad U(0)=U_0\in H.
			\end{equation}
			Suppose that $\cT:D(\cT)\subset H \arr H$ is a maximal monotone mapping,  $0\in \cT0$ and $B:H\arr H$ is locally Lipschitz, i.e. there exits $L(K)>0$ such that
			\begin{equation*}
				||B(U)-B(V)||_H \le L(K)||U-V||_H, \quad ||U||_H, ||V||_H \le K.
			\end{equation*}
			If $U_0\in D(\cT)$, $f\in W_1^1(0,t;H)$ for all $t>0$, then there exists $t_{max}\le \infty$ such that \eqref{AbstrIVP} has a unique strong solution $U$ on $(0,t_{max})$.\\
			If $U_0\in \overline{D(\cT)}$, $f\in L^1(0,t;H)$ for all $t>0$, then there exists $t_{max}\le \infty$ such that \eqref{AbstrIVP} has a unique generaized solution $U$ on $(0,t_{max})$.\\
			In both cases
			\begin{equation*}
				\lim_{t\arr t_{max}}||U(t)||_H=\infty \quad \mbox{provided} \quad t_{max}<\infty.
			\end{equation*}
		\end{theorem}
		First, we need to check that $\cT$ is a maximal monotone operator. Monotonicity is a direct consequence of Lemma \ref{lem:ASelfAdjont} and \eqref{DCont}. \\
		To prove $\cT$ is maximal as an operator from $H$ to $H$, we use Theorem 1.2 from \cite[Ch. 2]{Bar1976}. Thus, we need to prove that $Range(I+\cT)=H$. Let $z=(\Phi_z,\Psi_z)\in H_d\times H_v$. We need to find $y=(\Phi_y,\Psi_y)\in D(A)\times H_d=D(\cT)$ such that
		\begin{align*}
			& -\Psi_y + \Phi_y=\Phi_z,\\
			& A\Phi_y+\Psi_y + \Ga(\Psi_y) = \Psi_z,
		\end{align*}
		or, equivalently, find $\Psi_y\in H_d$ such that
		\begin{equation*}
			M(\Psi_y)=\frac 12 A\Psi_y + \frac 12 A\Psi_y+\Psi_y +\Ga(\Psi_y) = \Psi_z - A\Phi_z = \Theta_z
		\end{equation*}
		for an arbitrary $\Theta_z\in H_d'=D(A^{1/2})'$.
		Naturally, due to Lemma \ref{lem:ASelfAdjont} $A$ is a duality map between $H_d$ and $H_d'$, thus the operator $M$ is onto if and only if $\frac 12 A\Psi_y+\Psi_y +\Ga(\Psi_y)$ is maximal monotone as an operator from  $H_d$ to $H_d'$.
		According to Corollary 1.1 from \cite[Ch. 2]{Bar1976}, this operator is maximal monotone if $\frac 12 A$ is maximal monotone (it follows from Lemma \ref{lem:ASelfAdjont}) and $I+\Ga(\cdot)$ is monotone, bounded and hemicontinuous from  $H_d$ to $H_d'$. The last statement is evident for the identity map, now let's prove it for $\Ga$.\\
		Monotonicity is evident. Due to the continuity of the embedding $H^1(0,L_0)\subset C(0,L_0)$  in 1D every bounded  set $X$ in $H^1(0,L_0)$ is bounded in  $C(0,L_0)$ and thus due to \eqref{DCont} $\Ga(X)$ is bounded in $C(0,L_0)$ and, consequently, in $L^2(0,L_0)$. To prove hemicontinuity we take an arbitrary $\Phi=(\vph, \psi, \om, u,v,w) \in H_d$,  an  arbitrary $\Theta=(\theta_1,\theta_2, \theta_3, \theta_4,\theta_5, \theta_6)\in H_d$ and  consider
		\begin{equation*}
			(\Ga(\Psi_y+t\Phi),\Theta)=\int_0^{L_0}  \ga(\psi_y(x) + t\psi(x))\theta_2(x)dx,
		\end{equation*}
		where $\Psi_y=(\vph_y, \psi_y, \om_y, u_y,v_y,w_y)$. Since $\psi_y + t\psi\arr \psi_y, $ as $ t\arr 0$  in $H^1(0,L_0)$ and in  $C(0,L_0)$, we obtain that $\ga(\psi_y(x) + t\phi(x)) \arr \ga(\psi_y(x))$ as $t\arr 0$  for every $x\in [0,L_0]$, and has an integrable bound from above due to \eqref{DCont}. This implies $\ga(\psi_y(x) + t\phi(x)) \arr \ga(\psi_y(x))$ in $L^1(0,L_0)$ as $t\arr 0$ .  Since $\theta_2\in H^1(0,L_0)\subset L^\infty (0,L_0)$, then
		\begin{equation*}
			(\Ga(\Psi_y+t\Phi),\Theta) \arr (\Ga(\Psi_y),\Theta), \quad t\arr 0.
		\end{equation*}
		Hemicontinuity is proved.\\
		Further, we need to prove that $\cB$ is locally Lipschitz on $H$, that is, $F$ is locally Lipschitz from $H_d$ to $H_v$. The embedding $H^{1/2+\ep}(0,L)\subset C(0,L)$  and \eqref{NLip} imply
		\begin{equation}
			|F_j(\widetilde{\Phi}^j(x))-F_j(\widehat{\Phi}^j(x))|\le C(\max( ||\widetilde{\Phi}||_d, ||\widehat{\Phi}||_d)) ||\widetilde{\Phi}^j-\widehat{\Phi}^j||_{1}  \label{FLip}
		\end{equation}
		for all $x\in [0,L_0]$, if $j=1$ and for all $x\in [L_0,L]$, if $j=2$. This, in turn, gives us the estimate
		\begin{equation*}
			||F(\widetilde{\Phi})-F(\widehat{\Phi})||_v\le C(\max(||\widetilde{\Phi}||_d, ||\widehat{\Phi}||_d)) ||\widetilde{\Phi}-\widehat{\Phi}||_d.
		\end{equation*}
		Thus, all the assumptions of Theorem \ref{th:AbstrEx} are satisfied and existence of a local strong/generalized solution is proved.\\
		{\it Step 3. Energy inequality and  global solutions.}
		It can be verified by direct calculations, that strong solutions satisfy energy equality. Using the same arguments, as in proof of Proposition 1.3 \cite{ChuLa2007}, and \eqref{DCont} we can pass to the limit and prove \eqref{EE} for generalized solutions.  \\
		Let us assume that a local generalised solution exists on a maximal interval $(0, t_{max})$, $t_{max}<\infty$. Then \eqref{EE} implies $\cE(t_{max})\le \cE(0)$. Since due to \eqref{NBoundBelow}
		\begin{equation*}
			c_1||U(t)||_H \le \cE(t) \le c_2||U(t)||_H,
		\end{equation*}
		we have $||U(t_{max})||_H\le C||U_0||_H$. Thus, we arrive to a contradiction which implies $t_{max}=\infty$.
		
		{\it Step 4. Generalized solution is variational (weak).}  We formulate the following obvious estimate as a lemma for  future use.
		\begin{lemma}\label{lem:lip}
			Let \eqref{NLip} holds and $\widetilde{\Phi}$, $\widehat{\Phi}$ are two weak solutions to \eqref{AEq}-\eqref{AIC} with the initial conditions $(\widetilde{\Phi}_0, \widetilde{\Phi}_1)$ and  $(\widehat{\Phi}_0, \widehat{\Phi}_1)$ respectively. Then the following estimate is valid for all $x\in [0,L], \; t>0$ and $\epsilon\in[0,1/2)$:
			\begin{equation*}
				|F_j(\widetilde{\Phi}^j(x,t))-F_j(\widehat{\Phi}^j(x,t))|\le C(\max(||(\widetilde{\Phi}_0, \widetilde{\Phi}_1)||_H, ||(\widehat{\Phi}_0, \widehat{\Phi}_1)||_H)) ||\widetilde{\Phi}^j(\cdot,t)-\widehat{\Phi}^j(\cdot,t)||_{1-\epsilon}, \quad j=1,2.
			\end{equation*}
		\end{lemma}
		\begin{proof}
			The energy inequality and the embedding $H^{1/2+\ep}(0,L)\subset C(0,L)$  imply that for every weak solution $\Phi$
			\begin{equation*}
				\max_{t\in[0,T], x\in[0,L]}	|\Phi(x,t)| \le C(||\Phi_0||_d, ||\Phi_1||_v).
			\end{equation*}
			Thus, using \eqref{NLip}  and  \eqref{FLip}, we prove the lemma.
		\end{proof}	
		Evidently, \eqref{VEq} is valid for strong solutions. We can find a sequence of strong solutions $\Phi^{(n)}$, which converges to a generalized solution $\Phi$ strongly in $C(0,T; H_d)$, and $\Phi_t^{(n)}$ converges to $\Phi_t$ strongly in $C(0,T; H_v)$. Using Lemma \ref{lem:lip}, we can easily pass to the  limit  in nonlinear feedback term in \eqref{VEq}. Since the test function  $B\in L^\infty(0,T;H_d)\subset L^\infty((0,T)\times (0,L))$, we can use the same arguments as in the proof of Proposition 1.6 \cite{ChuLa2007} to pass to the limit in the nonlinear dissipation term.  Namely, we can extract from $\Phi_t^{(n)}$  a subsequence that  converges to $\Phi_t$  almost everywhere and prove that it converges to $\Phi_t$ strongly in $L^1((0,T)\times (0,L))$.
	\end{proof}
	\begin{remark}
		In space dimension greater then one we do not have the embedding $H^1(\Omega)\subset C(\Omega)$, therefore, we need to assume polynomial growth of the derivative of the nonlinearity to obtain estimates similar to Lemma \ref{lem:lip}.
	\end{remark}
	\section{Existence of attractors.}
	In this section we study long-time behaviour of solutions to problem \eqref{AEq}-\eqref{AIC} in the framework of dynamical systems theory. From Theorem \ref{th:WeakWP} we have
	\begin{corollary}\label{DSGen}
		Let, additionally to conditions of Theorem \ref{th:WeakWP}, $P(x,t)=P(x)$. Then \eqref{AEq}-\eqref{AIC} generates a dynamical system $(H, S_t)$ by the formula
		\begin{equation*}
			S_t(\Phi_0,\Phi_1)=(\Phi(t),\Phi_t(t)),
		\end{equation*}
		where $\Phi(t)$ is the weak solution to \eqref{AEq}-\eqref{AIC} with initial data $(\Phi_0,\Phi_1)$.
	\end{corollary}	
	To establish the  existence of the attractor for this dynamical system  we use Theorem \ref{abs} below, thus we need to prove the gradientness and  the asymptotic smoothness as well as the boundedness of the set of stationary points.
	\subsection{Gradient structure}
	In this subsection we prove that the dynamical system generated by \eqref{AEq}-\eqref{AIC} possesses a specific structure, namely, is gradient under some additional conditions on the nonlinearities.
	\begin{definition}[\cite{Chueshov,CFR,CL}]\label{de:grad}
		Let $Y\subseteq X$ be a positively invariant set of  $(X,S_t)$.
		\begin{itemize}
			\item a continuous functional  $L(y)$ defined on  $Y$ is said to be a \emph{Lyapunov function} of the dynamical system  $(X,S_t)$ on the set $Y$, if
			a function  $t\mapsto L(S_ty)$ is non-increasing for any $y\in Y$.
			\item
			the Lyapunov function  $L(y)$ is said to be \emph{strict} on $Y$,
			if the equality  $L(S_{t}y)=L(y)$  \emph{for all} $t>0$ implies  $S_{t}y=y$ for all $t>0$;
			\item
			A dynamical system $(X,S_t)$ is said to be \emph{gradient}, if it possesses a strict Lyapunov function on the whole phase space $X$.
		\end{itemize}
	\end{definition}
	The following result holds true.
	\begin{theorem}\label{th:grad}
		Let, additionally to the assumptions of Corollary \ref{DSGen}, the following conditions hold
		\begin{align*}
			& f_1=g_1=0, \quad h_1(\vph, \psi,\om)=h_1(\psi), \label{NZero}\tag{N3} \\
			& f_2,\; g_2,\; h_2 \in C^1(\R^3), \label{NAddSmooth}\tag{N4} \\
			& \ga(s)s > 0 \quad \mbox{ for all } s\neq 0. \label{DNonZero} \tag{D2}
		\end{align*}
		Then the dynamical system $(H, S_t)$ is gradient.
	\end{theorem}
	\begin{proof}
		We use as a Lyapunov function
		\begin{equation}\label{lap}
			L(\Phi(t))=L(t) =\frac 12\left( ||R^{1/2}\Phi_t(t)||^2+ ||A^{1/2}\Phi(t)||^2\right) +
			\intl_0^L  \cF(\Phi(x,t))dx + (P,\Phi(t)).
		\end{equation}
		Energy inequality \eqref{EE} implies that $L(t)$ is non-increasing. The equality $L(t)=L(0)$ together with \eqref{DNonZero} imply that $\psi_t(t)\equiv 0$ on $[0,T]$. We need to prove that  $\Phi(t)\equiv const$, which is equivalent to $\Phi(t+h)-\Phi(t)=0$ for every $h>0 $. In what follows we use the notation $\Phi(t+h)-\Phi(t)=\oPhi(t)=(\ovph, \overline{\psi},\oom,\ou,\ov,\ow)(t)$ .\\
		{\it Step 1.} Let us  prove that $\oPhi^1\equiv 0$. In this step we use the distribution theory (see, e.g., \cite{Kan2004}) because some functions involved in computations are of  too low smoothness. Let us set the test function $B=(B^1,0)=(\beta^1,\gamma^1,\delta^1, 0,0,0)$. Then $\oPhi(t)$ satisfies
		\begin{multline*}
			-\intl_0^T  (R_1 \oPhi^1_t, B_t)(t)dt - (R_1(\Phi^1_t(h)-\Phi_1^1), B^1(0)) +\\
			\intl_0^T \left[\frac 1{k_1}(Q_1(\oPhi^1),Q_1(B^1))(t) dt+ \frac 1{\si_1}(N_1(\oPhi^1), N_1(B^1))(t) \right] + \\
			\intl_0^T (h_1(\psi(t+h))-h_1(\psi(t)), \gamma^1(t))dt =0.
		\end{multline*}
		The last term equals to zero due to \eqref{NZero}  and $\psi(t)\equiv const$.\\
		Setting  in turn $B=(0,\gamma^1,0, 0,0,0)$, $B=(0,0,\de^1, 0,0,0)$, $B=(\beta^1, 0,0, 0,0,0)$ we obtain
		\begin{align}
			& \ovph_x+l\oom=0 \qquad&\mbox{ almost everywhere on } (0,L_0)\times (0,T),  \label{GrSpace} \\
			& \rho_1\oom_{tt} - l\si_1(\oom_x-l\ovph)_x =0  \qquad&\mbox{ almost everywhere on } (0,L_0)\times (0,T), \label{GrOmega} \\
			& \rho_1\ovph_{tt} - \si_1(\oom_x-l\ovph)=0  \qquad&\mbox{ in the sense of distributions on } (0,L_0)\times (0,T). \label{GrPhi}
		\end{align}
		These equalities imply
		\begin{equation}\label{GrTime}
			\ovph_{ttx}=0, \quad \oom_{tt}=0 \quad \mbox{ in the sense of distributions}.
		\end{equation}
		Similar to regular functions,  if partial derivative of a distribution equals to zero, then the distribution "does not depends" on the corresponding variable (see \cite[Ch. 7]{Kan2004}, Example 2.) That is,
		\begin{equation*}
			\oom_t=c_1(x)\times 1(t) \qquad \mbox{ in the sense of  distributions} \\
		\end{equation*}
		However, Theorem \ref{th:WeakWP} implies that  $\oom_t$ is a regular distribution, thus, we can treat the equality above as the equality almost everywhere. Furthermore,
		\begin{equation*}
			\oom(x,t)=\oom(x,0)+\int_0^t  c_1(x)d\tau = \oom(x,0) +tc_1(x).
		\end{equation*}
		Since $||\oom(\cdot,t)||\le C$ for all $t\in\R_+$, $c_1(x)$ must be zero. Thus,
		\begin{equation}\label{GrOmegaConst}
			\oom(x,t)=c_2(x),
		\end{equation}
		which together with \eqref{GrSpace} implies
		\begin{align*}
			& \ovph_x=-lc_2(x), \\
			& \ovph(x,t)= \ovph(0,t) - l\intl_0^x  c_2(y)dy = c_3(x), \\
			& \ovph_{tt} =0.
		\end{align*}
		The last equality together with \eqref{GrSpace}, \eqref{GrPhi} and boundary conditions \eqref{ABC} give us that $\ovph, \oom$ are solutions to the following Cauchy problem (with respect to $x$):
		\begin{align*}
			&\oom_x = l\ovph,\\
			&\ovph_x = -l\oom, \\
			&\oom(0,t)=\ovph(0,t)=0.
		\end{align*}
		Consequently, $\oom\equiv\ovph\equiv 0$.\\
		{\it Step 2.} Let us prove, that $u\equiv v\equiv w\equiv 0$. Due to \eqref{NAddSmooth}, we can use the Taylor expansion of the difference $F^2(\Phi^2(t+h))-F^2(\Phi^2(t))$  and thus $(\ou, \ov, \ow)$ satisfies
		on $(0,T)\times (L_0,L)$
		\begin{align}
			& \rho_2\ou_{tt}-k_2\ou_{xx} +g_u(\pd_x \oPhi^2, \oPhi^2) + \g f_2(\zeta_{1,h}(x,t))\cdot\oPhi^2=0, \label{WE1} \\
			& \beta_2\ov_{tt} -\la_2 \ov_{xx} + g_v(\pd_x \oPhi^2, \oPhi^2) + \g h_2(\zeta_{2,h}(x,t))\cdot\oPhi^2=0,\\
			& \rho_2\ow_{tt}- \si_2 \ow_{xx} + g_w(\pd_x \oPhi^2, \oPhi^2) + \g g_2(\zeta_{3,h}(x,t))\cdot\oPhi^2=0 \\
			& \ou(L_0,t)=\ov(L_0,t)=\ow(L_0,t)=0,\\
			& \ou(L,t)=\ov(L,t)=\ow(L,t)=0, \\
			& {k_2(\ou_x+\ov+l\ow)(L_0,t)=0}, \label{WEBC1}\\
			&{\ov_x(L_0,t)=0, \qquad \sigma_2(\ow_x-l\ou)(L_0,t)=0}, \label{WEBC2}\\
			& \oPhi^2(x,0)=\Phi^2(x,h)-\Phi^2_0, \quad \oPhi^2_t(x,0)=\Phi^2_t(x,h)-\Phi^2_1, \label{WEIC}
		\end{align}
		where $g_u, g_v, g_w$ are linear combinations of $u_x,v_x,w_x, u,v,w$ with the constant coefficients, $\zeta_{j,h}(x,t)$ are 3D vector functions which components lie between $u(x, t+h)$ and $u(x,t)$, $v(x,t+h)$ and $v(x,t)$, $w(x,t+h)$ and $w(x,t)$ respectively. Thus, we have a system of linear equations on $(L_0,L)$ with overdetermined boundary conditions. $L^2$-regularity of $u_x,v_x, w_x$ on the boundary for solutions to a linear wave equation  was established in \cite{LaTri1983}, thus, boundary conditions \eqref{WEBC1}-\eqref{WEBC2} make sense.\\
		It is easy to generalize the observability inequality \cite[Th. 8.1]{TriYao2002} for the case of the system of the wave equations.
		\begin{theorem}[\cite{TriYao2002} ]
			For the solution to problem \eqref{WE1}-\eqref{WEIC} the following estimate holds:
			\begin{equation*}
				\int_0^T  [|\ou_x|^2+|\ov_x|^2+|\ow_x|^2](L_0,t) dt\ge C(E(0)+E(T)),
			\end{equation*}
			where
			\begin{equation*}
				E(t)=\frac 12 \left(||\ou_t(t)||^2+ ||\ov_t(t)||^2 +||\ow_t(t)||^2 + ||\ou_x(t)||^2 +||\ov_x(t)||^2+||\ow_x(t)||^2 \right).
			\end{equation*}
		\end{theorem}
		Therefore, if conditions \eqref{WEBC1}, \eqref{WEBC2} hold true, then $\ou=\ov=\ow=0$.
		The theorem is proved.
	\end{proof}
	\subsection{Asymptotic smoothness.}
	\begin{definition}[\cite{Chueshov,CFR,CL}]
		A dynamical system $(X,S_t)$ is said to be asymptotically smooth
		if for any  closed bounded set $B\subset X$ that is positively invariant ($S_tB\subseteq B$)
		one can find a compact set $\cK=\cK(B)$ which uniformly attracts $B$, i.~e.
		$\sup\{{\rm dist}_X(S_ty,\cK):\ y\in B\}\to 0$ as $t\to\infty$.
	\end{definition}
	In order to prove  the asymptotical smoothness of the system
	considered we rely on the  compactness criterion due to
	\cite{Khanmamedov}, which is recalled below in an abstract version
	formulated in \cite{CL}.
	\begin{theorem}{\cite{CL}} \label{theoremCL} Let $(S_t, H)$ be a dynamical system on a complete metric
		space $H$ endowed with a metric $d$. Assume that for any bounded positively invariant
		set $B$ in $H$ and for any $\varepsilon>0$ there exists $T = T (\varepsilon, B)$ such that
		\begin{equation}
			\label{te}
			d(S_T y_1, S_T y_2) \le \varepsilon+ \Psi_{\varepsilon,B,T} (y_1, y_2), y_i \in B ,
		\end{equation}
		where $\Psi_{\varepsilon,B,T} (y_1, y_2)$ is a function defined on $B \times B$ such that
		\[
		\liminf\limits_{m\to\infty}\liminf\limits_{n\to\infty}\Psi_{\varepsilon,B,T} (y_n, y_m) = 0
		\]
		for every sequence ${y_n} \in B$. Then $(S_t, H)$ is an asymptotically smooth dynamical
		system.
	\end{theorem}
	
	To formulate the result on the asymptotic smoothness of the system considered we need the following lemma.
	\begin{lemma}
		\label{lem:GammaEst}
		Let assumptions \eqref{DCont} hold. Let moreover,
		there exists a positive constant $M$ such that
		\begin{equation}
			\frac{\gamma(s_1)-\gamma(s_2)}{s_1-s_2}\le M, \quad s_1, s_2\in \R, \,\,s_1\ne s_2.\tag{D3}\label{GammaLip1}
		\end{equation}
		Then for any $\varepsilon>0$ there exists $C_\varepsilon>0$ such that
		\begin{equation}\label{GammaEst}
			\left|\int\limits_0^{L_0}(\ga( \xi_1)-\ga( \xi_2)) \zeta dx\right|\le \varepsilon \|\zeta\|^2+C_\varepsilon \int\limits_0^{L_0}(\ga( \xi_1)-\ga( \xi_2)) (\xi_1-\xi_2)dx
		\end{equation}
		for any $\xi_1, \xi_2, \zeta\in L^2(0,L_0)$.
	\end{lemma}
	The proof is similar to  that given in \cite[Th.5.5]{CL}).
	\begin{theorem}
		\label{th:AsSmooth}
		Let assumptions of Theorem \ref{th:WeakWP}, \eqref{GammaLip1}, and
		\begin{equation}
			m\le \frac{\gamma(s_1)-\gamma(s_2)}{s_1-s_2}, \quad s_1, s_2\in \R, \,\,s_1\ne s_2.\tag{D4}\label{GammaLip2}
		\end{equation}
		with $m>0$	 hold. Let, moreover,
		\begin{gather}
			k_1=\sigma_1 \label{c1}\\
			\frac{\rho_1}{k_1}=\frac{\beta_1}{\lambda_1}. \label{c2}
		\end{gather}
		Then the dynamical system  $(H, S_t)$ generated by problem \eqref{Eq1}--\eqref{TC4} is asymptotically smooth.
	\end{theorem}
	\begin{proof}
		In this proof we perform all the calculations for strong solutions and then pass to the limit in the final estimate to justify it for weak solutions.
		Let us consider strong solutions $\hat U(t)=(\hat\Phi(t), \hat\Phi_t(t))$ and $\tilde U(t)=(\tilde\Phi(t), \tilde\Phi_t(t))$  to problem \eqref{Eq1}--\eqref{TC4} with initial conditions $\hat U_0=(\hat \Phi_0, \hat \Phi_1)$ and $\tilde U_0=(\tilde \Phi_0, \tilde \Phi_1)$ lying in a ball, i.e.  there exists $R>0$ such that
		\begin{equation}
			\label{inboun}
			\|\tilde U_0\|_H+\|\hat U_0\|_H\le R.
		\end{equation}
		Denote $U(t)=\tilde U(t)-\hat U(t)$ and  $U_0=\tilde U_0-\hat U_0$. Obviously, $U(t)$ is a weak solution to the problem
		\begin{align}
			& \rho_1\vph_{tt}-k_1(\vph_x+\psi+l\om)_x - l\si_1(\om_x-l\vph) +f_1(\tilde\vph, \tilde\psi, \tilde\om)-f_1(\hat\vph, \hat\psi, \hat\om)=0\label{eq1}\\
			& \be_1\psi_{tt} -\la_1 \psi_{xx} +k_1(\vph_x+\psi+l\om) +\ga( \tilde\psi_t)-\ga( \hat\psi_t) +h_1(\tilde\vph, \tilde\psi, \tilde\om)-h_1(\hat\vph, \hat\psi, \hat\om)=0\label{eq2}\\
			& \rho_1\om_{tt}- \si_1(\om_x-l\vph) _x+lk_1(\vph_x+\psi+l\om)+g_1(\tilde\vph, \tilde\psi, \tilde\om)-g_1(\hat\vph, \hat\psi, \hat\om)=0 \label{eq3}\\		
			& \rho_2u_{tt}-k_2(u_x+v+lw)_x - l\si_2(w_x-lu) +f_2(\tilde u, \tilde v, \tilde w)-f_2(\hat u, \hat v, \hat w)=0 \label{eq4}\\
			& \be_2v_{tt} -\la_2 v_{xx} +k_2(u_x+v+lw) +h_2(\tilde u, \tilde v, \tilde w)-h_2(\hat u, \hat v, \hat w)=0,\label{eq5}\\
			& \rho_2w_{tt}- \si_2(w_x-lu)_x+lk_2(u_x+v+lw)+g_2(\tilde u, \tilde v, \tilde w)-g_2(\hat u, \hat v, \hat w)=0 \label{eq6}
		\end{align}
		with boundary conditions \eqref{BC}, \eqref{TC1}--\eqref{TC4} and the initial conditions $U(0)=\tilde U_0-\hat U_0$.
		It is easy to see by the energy argument that
		\begin{equation}
			\label{En1}
			E(U(T))+ \int\limits_t^T \int\limits_0^{L_0}(\ga( \tilde\psi_s)-\ga( \hat\psi_s)) \psi_s dx ds=E(U(t))+ \int\limits_t^T  H(\hat U(s),\tilde U(s))  ds,
		\end{equation}
		where
		\begin{multline}
			\label{H}
			H(\hat U(t),\tilde U(t))=\int\limits_0^{L_0}(f_1(\hat\vph, \hat\psi, \hat\om)-f_1(\tilde\vph, \tilde\psi, \tilde\om))\vph_t dx+\int\limits_0^{L_0}(h_1(\hat\vph, \hat\psi, \hat\om)-h_1(\tilde\vph, \tilde\psi, \tilde\om))\psi_t dx\\+
			\int\limits_0^{L_0}(g_1(\hat\vph, \hat\psi, \hat\om)-g_1(\tilde\vph, \tilde\psi, \tilde\om))\omega_t dx+
			\int\limits_{L_0}^L(f_2(\hat u, \hat v, \hat w)-f_2(\tilde u, \tilde v, \tilde w)) u_t dx\\+\int\limits_{L_0}^L(h_2(\hat u, \hat v, \hat w)-h_2(\tilde u, \tilde v, \tilde w))v_t dx+
			\int\limits_{L_0}^L(g_2(\hat u, \hat v, \hat w)-g_2(\tilde u, \tilde v, \tilde w))w_t dx,
		\end{multline}
		and
		\begin{equation}
			\label{E}
			E(t)=E_1(t)+E_2(t),
		\end{equation}
		here
		\begin{multline}
			\label{E1}
			E_1(t)=\rho_1\int\limits_0^{L_0} \omega_t^2dx dt+\rho_1 \int\limits_0^{L_0}\vph_{t}^2 dx dt+\beta_1 \int\limits_0^{L_0}\psi_{t}^2 dx+\sigma_1\int\limits_0^{L_0} (\omega_x-l\vph)^2 dx+\\+k_1\int\limits_0^{L_0} (\vph_x+\psi+l\om)^2 dx+\lambda_1\int\limits_0^{L_0} \psi_x^2 dx
		\end{multline}
		and
		\begin{multline}
			\label{E2}
			E_2(t)=\rho_2\int\limits_0^{L_0} w_t^2dx dt+\rho_2 \int\limits_0^{L_0}u_{t}^2 dx dt+\beta_2 \int\limits_0^{L_0}v_{t}^2 dx+\sigma_2\int\limits_0^{L_0} (w_x-lu)^2 dx+\\+k_2\int\limits_0^{L_0} (u_x+v+lw)^2 dx+\lambda_2\int\limits_0^{L_0} v_x^2 dx.
		\end{multline}
		Integrating in \eqref{En1}  over the interval $(0,T)$ we come to
		\begin{equation}
			\label{En2}
			TE(U(T))+\int\limits_0^T  \int\limits_t^T \int\limits_0^{L_0} (\ga( \tilde\psi_s)-\ga( \hat\psi_s))\psi_s dx ds dt\\=\int\limits_0^T  E(U(t)) dt+ \int\limits_0^T \int\limits_t^T  H(\hat U(s),\tilde U(s))  ds dt.
		\end{equation}
		Now we  estimate the first term in the right-hand side of \eqref{En2}. In what follows we present formal estimates which can be performed on strong solutions.\\
		{\it Step 1.} We multiply equation \eqref{eq3} by $\omega$ and $x\cdot\omega_x$ and sum up the results. After integration by parts with respect to $t$ we obtain
		\begin{multline*}
			\rho_1 \int\limits_0^T \int\limits_0^{L_0} \omega_t x \omega_{tx}  dx dt +\rho_1 \int\limits_0^T \int\limits_0^{L_0} \omega_t^2  dx dt+\sigma_1 \int\limits_0^T \int\limits_0^{L_0} (\omega_x-l\vph)_x x \omega_{x}  dx dt\\
			+\sigma_1 \int\limits_0^T \int\limits_0^{L_0} (\omega_x-l\vph)_x \omega  dx dt-k_1 l \int\limits_0^T \int\limits_0^{L_0}(\vph_x+\psi+l\omega)x\omega_x dx dt\\
			-k_1 l \int\limits_0^T \int\limits_0^{L_0}(\vph_x+\psi+l\omega)\omega dx dt- \int\limits_0^T \int\limits_0^{L_0}(g_1(\tilde\vph,\tilde\psi,\tilde\omega)-g_1(\hat\vph,\hat\psi,\hat\omega))(x\omega_x+\omega) dx dt\\
			=\rho_1  \int\limits_0^{L_0} \omega_t(x,T) x \omega_{x}(x,T)  dx+\rho_1 \int\limits_0^{L_0} \omega_t(x,T) \omega(x,T) dx
		\end{multline*}
		\begin{equation}\label{o1}
			-\rho_1  \int\limits_0^{L_0} \omega_t(x,0) x \omega_{x}(x,0)  dx-\rho_1 \int\limits_0^{L_0} \omega_t(x,0) \omega(x,0) dx.
		\end{equation}
		Integrating by parts with respect to $x$ we get
		\begin{equation}
			\label{o2}
			\rho_1 \int\limits_0^T \int\limits_0^{L_0} \omega_t x \omega_{tx}  dx dt =-\frac{\rho_1}{2} \int\limits_0^T \int\limits_0^{L_0} \omega_t^2  dx dt+\frac{\rho_1L_0}{2}\int\limits_0^T\omega_t^2(L_0,t)dt
		\end{equation}
		and
		\begin{multline}	\label{o3}
			\sigma_1 \int\limits_0^T \int\limits_0^{L_0} (\omega_x-l\vph)_x x \omega_{x} dx dt -k_1 l \int\limits_0^T \int\limits_0^{L_0}(\vph_x+\psi+l\omega)x\omega_x dx dt\\
			= \sigma_1 \int\limits_0^T \int\limits_0^{L_0} (\omega_x-l\vph)_x x (\omega_x-l\vph) dx dt+\sigma_1 l\int\limits_0^T \int\limits_0^{L_0} (\omega_x-l\vph)_x x \vph dx dt\\-k_1 l \int\limits_0^T \int\limits_0^{L_0}(\vph_x+\psi+l\omega)x\omega_x dx dt=
			-	\frac{\sigma_1}{2} \int\limits_0^T \int\limits_0^{L_0} (\omega_x-l\vph)^2 dx dt\\
			+\frac{\sigma_1 L_0}{2} \int\limits_0^T (\omega_x-l\vph)^2(L_0,t) dt-\sigma_1 l\int\limits_0^T \int\limits_0^{L_0} (\omega_x-l\vph)  \vph dx dt\\
			-2\sigma_1 l\int\limits_0^T \int\limits_0^{L_0} (\omega_x-l\vph)x(\vph_x +\psi+l\omega)dx dt+\sigma_1 l\int\limits_0^T \int\limits_0^{L_0} (\omega_x-l\vph)x(\psi+l\omega)dx dt\\
			-\sigma_1 lL_0\int\limits_0^T  (\omega_x-l\vph) (L_0,t) \vph (L_0,t) dt-k_1 l^2 \int\limits_0^T \int\limits_0^{L_0}(\vph_x+\psi+l\omega)x\vph dx dt.
		\end{multline}
		Analogously,
		\begin{multline}	\label{o4}
			\sigma_1 \int\limits_0^T \int\limits_0^{L_0} (\omega_x-l\vph)_x \omega  dx dt=-\sigma_1 \int\limits_0^T \int\limits_0^{L_0} (\omega_x-l\vph)^2  dx dt\\
			+\sigma_1 \int\limits_0^T  (\omega_x-l\vph)(L_0,t) \omega(L_0,t) dt-l\sigma_1 \int\limits_0^T \int\limits_0^{L_0} (\omega_x-l\vph)\vph  dx dt.
		\end{multline}
		It follows from Lemma \ref{lem:lip}, energy relation \eqref{EE}, and property \eqref{NBoundBelow}  that
		\begin{equation}
			\label{nn}
			\int\limits_0^T \int\limits_0^{L_0}|g_1(\tilde\vph,\tilde\psi,\tilde\omega)-g_1(\hat\vph,\hat\psi,\hat\omega)|^2 dx dt \le C(R,T)\max\limits_{t\in[0,T]} \|\Phi(\cdot,t)\|_{H^{1-\epsilon}}^2,\, 0<\epsilon<1/2.
		\end{equation}
		Therefore, for every $\varepsilon>0$
		\begin{equation}	\label{o5}
			\left|\int\limits_0^T \int\limits_0^{L_0}(g_1(\tilde\vph,\tilde\psi,\tilde\omega)-g_1(\hat\vph,\hat\psi,\hat\omega))(x\omega_x+\omega) dx dt\right|\\
			\le \varepsilon\int\limits_0^T\|\omega_x-l\vph\|^2 dt+C(\varepsilon,R,T) lot,
		\end{equation}
		where we use  the notation
		\begin{multline}
			\label{lot}
			lot=\max\limits_{t\in[0,T]} (\| \vph(\cdot,t)\|_{H^{1-\epsilon}}^2+\| \psi(\cdot,t)\|_{H^{1-\epsilon}}^2+ \|\omega(\cdot,t)\|_{H^{1-\epsilon}}^2\\
			+\| u(\cdot,t)\|_{H^{1-\epsilon}}^2+ \|v(\cdot,t)\|_{H^{1-\epsilon}}^2+\|w(\cdot,t)\|_{H^{1-\epsilon}}^2),\quad 	0<\epsilon<1/2.
		\end{multline}
		Similar estimates hold for  nonlinearities $g_2$, $f_i$, $h_i$, $i=1,2$.\\
		We note that for any $\eta\in H^1(0,L_0)$ (or analogously  $\eta\in H^1(L_0, L)$)
		\begin{equation}
			\label{lote}
			\eta(L_0)\le \sup\limits_{(0,L_0)}|\eta|\le C \| \eta\|_{H^{1-\epsilon}}, \quad 	0<\epsilon<1/2.		
		\end{equation}
		Since due to \eqref{c1}
		\begin{multline*}
			2\sigma_1 l\left|\int\limits_0^T \int\limits_0^{L_0} (\omega_x-l\vph)x(\vph_x +\psi+l\omega)dx dt\right|\\\le 	\frac{\sigma_1}{16} \int\limits_0^T \int\limits_0^{L_0} (\omega_x-l\vph)^2 dx dt+16k_1 l^2L_0^2\int\limits_0^T \int\limits_0^{L_0} (\vph_x +\psi+l\omega)^2dx dt,
		\end{multline*}
		the following estimate can be obtained from \eqref{o1}-- \eqref{o5}
		\begin{multline}	\label{ol1}
			\frac{\rho_1}{2} \int\limits_0^T \int\limits_0^{L_0} \omega_t^2dx dt +\frac{\rho_1L_0}{2}\int\limits_0^T\omega_t^2(L_0,t)dt+\frac{13\sigma_1L_0}{8} \int\limits_0^T  (\omega_x-l\vph)^2(L_0,t) dt\\
			\le \frac{13\sigma_1}{8} \int\limits_0^T \int\limits_0^{L_0} (\omega_x-l\vph)^2  dx dt+17k_1 l^2L_0^2 \int\limits_0^T \int\limits_0^{L_0}(\vph_x+\psi+l\omega)^2 dx dt+C(R,T)lot\\+C(E(0)+E(T)),
		\end{multline}
		where $C>0$.\\
		{\it Step 2.} Multiplying  equation \eqref{eq3} by $\omega$ and $(x-L_0)\cdot\omega_x$ and arguing as above we come to the estimate
		\begin{multline}
			\label{ol2}
			\frac{\rho_1}{2} \int\limits_0^T \int\limits_0^{L_0} \omega_t^2dx dt+\frac{13\sigma_1L_0}{8} \int\limits_0^T  (\omega_x-l\vph)^2(0,t) dt
			\le \frac{13\sigma_1}{8} \int\limits_0^T \int\limits_0^{L_0} (\omega_x-l\vph)^2  dx dt\\+17k_1 l^2L_0^2 \int\limits_0^T \int\limits_0^{L_0}(\vph_x+\psi+l\omega)^2 dx dt+C(R,T)lot+C(E(0)+E(T)).
		\end{multline}
		Summing up estimates \eqref{ol1} and \eqref{ol2} and multiplying the result by $\frac{1}{2}$ we get
		
		\begin{multline}
			\label{ol3}
			\frac{\rho_1}{2}\int\limits_0^T \int\limits_0^{L_0} \omega_t^2dx dt +\frac{\rho_1L_0}{4}\int\limits_0^T\omega_t^2(L_0,t)dt+\frac{3\sigma_1L_0}{16} \int\limits_0^T  (\omega_x-l\vph)^2(L_0,t) dt\\+\frac{3\sigma_1L_0}{16} \int\limits_0^T  (\omega_x-l\vph)^2(0,t) dt\le \frac{13\sigma_1}{8} \int\limits_0^T \int\limits_0^{L_0} (\omega_x-l\vph)^2  dx dt\\+17k_1 l^2L_0^2 \int\limits_0^T \int\limits_0^{L_0}(\vph_x+\psi+l\omega)^2 dx dt+C(R,T)lot+C(E(0)+E(T)).
		\end{multline}
		{\it Step 3.}  Next we multiply equation \eqref{eq1} by $-\frac{1}{l}(\omega_x-l\vph)$, equation \eqref{eq3} by $\frac{1}{l}\vph_x$, summing up the results and integrating by parts with respect to $t$ we arrive at
		\begin{multline*}	
			\frac{\rho_1}{l}\int\limits_0^T \int\limits_0^{L_0}\vph_{t}(\omega_{tx}-l\vph_t)dx dt+\frac{k_1}{l}\int\limits_0^T \int\limits_0^{L_0}(\vph_x+\psi+l\om)_x (\omega_x-l\vph) dx dt\\ + \si_1\int\limits_0^T \int\limits_0^{L_0}(\om_x-l\vph)^2 dx dt-
			\frac{1}{l}\int\limits_0^T \int\limits_0^{L_0}(f_1(\tilde\vph,\tilde\psi,\tilde\omega)- f_1(\hat\vph,\hat\psi,\hat\omega))(\omega_x-l\vph) dx dt\\
			+\frac{\rho_1}{l}\int\limits_0^T \int\limits_0^{L_0}\om_{t}\vph_{tx} dx dt+ \frac{\si_1}{l}\int\limits_0^T \int\limits_0^{L_0}(\om_x-l\vph)_x\vph_{x} dx dt\\
			-k_1\int\limits_0^T \int\limits_0^{L_0}(\vph_x+\psi+l\om)\vph_{x} dx dt-\int\limits_0^T \int\limits_0^{L_0}(g_1(\tilde\vph,\tilde\psi,\tilde\omega)-g_1(\hat\vph,\hat\psi,\hat\omega))\vph_{x} dx dt\\
			=\frac{\rho_1}{l} \int\limits_0^{L_0}\vph_{t}(x,T)(\omega_{x}-l\vph)(x,T)dx-\frac{\rho_1}{l} \int\limits_0^{L_0}\vph_{t}(x,0)(\omega_{x}-l\vph)(x,0)dx
		\end{multline*}
		\begin{equation}	\label{o6}
			\qquad\qquad\qquad\qquad	+\frac{\rho_1}{l} \int\limits_0^{L_0}\om_{t}(x,T)\vph_{x}(x,T) dx-\frac{\rho_1}{l} \int\limits_0^{L_0}\om_{t}(x,0)\vph_{x}(x,0) dx.
		\end{equation}
		Integrating by parts with respect to $x$ we obtain
		\begin{equation*}
			\left|\frac{\rho_1}{l}\int\limits_0^T \int\limits_0^{L_0}\vph_{t}\omega_{tx}dx dt+\frac{\rho_1}{l}\int\limits_0^T \int\limits_0^{L_0}\om_{t}\vph_{tx} dx dt\right|
			=\left|\frac{\rho_1}{l}\int\limits_0^T \vph_{t}(L_0,t)\omega_{t}(L_0,t) dt\right|
		\end{equation*}
		\begin{equation}
			\label{o7}
			\qquad\qquad	\qquad\qquad	\qquad\qquad	\qquad\qquad\le	\frac{\rho_1L_0}{8}\int\limits_0^T\omega_{t}^2(L_0,t) dt+\frac{2\rho_1}{l^2L_0}\int\limits_0^T \vph_{t}^2(L_0,t) dt.
		\end{equation}
		Taking into account \eqref{c1} we get
		\begin{multline}
			\label{o8}
			\frac{k_1}{l}\int\limits_0^T \int\limits_0^{L_0}(\vph_x+\psi+l\om)_x (\omega_x-l\vph) dx dt+\frac{\si_1}{l}\int\limits_0^T \int\limits_0^{L_0}(\om_x-l\vph)_x\vph_{x} dx dt\\
			=\frac{k_1}{l}\int\limits_0^T (\vph_x+\psi+l\om)(L_0,t) (\omega_x-l\vph)(L_0,t)dt
			-\frac{k_1}{l}\int\limits_0^T (\vph_x+\psi+l\om)(0,t)(\omega_x-l\vph)(0,t)dt\\
			+\frac{k_1}{l}\int\limits_0^T \int\limits_0^{L_0}\psi_x (\omega_x-l\vph) dx dt+\si_1\int\limits_0^T \int\limits_0^{L_0}(\omega_x-l\vph)^2 dx dt+\si_1 l\int\limits_0^T \int\limits_0^{L_0}(\omega_x-l\vph)\vph dx dt.
		\end{multline}
		Using the estimates
		\begin{multline*}
			\left|\frac{k_1}{l}\int\limits_0^T (\vph_x+\psi+l\om)(L_0,t) (\omega_x-l\vph)(L_0,t)dt\right|\\\le
			\frac{4k_1}{l^2L_0}\int\limits_0^T(\vph_x+\psi+l\om)^2(L_0,t)dt+
			\frac{\si_1L_0}{16}\int\limits_0^T(\omega_x-l\vph)^2(L_0,t)dt,
		\end{multline*}
		\begin{equation*}
			\left|\frac{k_1}{l}\int\limits_0^T \int\limits_0^{L_0}\psi_x (\omega_x-l\vph) dx dt\right|\le
			\frac{4k_1}{l^2}\int\limits_0^T \int\limits_0^{L_0}\psi_x^2  dx dt+	\frac{\si_1}{16}\int\limits_0^T \int\limits_0^{L_0} (\omega_x-l\vph)^2 dx dt	
		\end{equation*}
		and \eqref{o6}--\eqref{o8}
		we infer
		
		\begin{multline}
			\label{ol4}
			\frac{15\si_1}{8}\int\limits_0^T \int\limits_0^{L_0} (\omega_x-l\vph)^2 dx dt\le  \rho_1 \int\limits_0^T \int\limits_0^{L_0} \vph_t^2 dx dt+ 2k_1\int\limits_0^T \int\limits_0^{L_0} (\vph_{x}+\psi+l\omega)^2dx dt\\
			+\frac{4k_1}{l^2} \int\limits_0^T \int\limits_0^{L_0} \psi_x^2 dx dt+\frac{4k_1}{l^2L_0}\int\limits_0^T  (\vph_{x}+\psi+l\omega)^2(L_0,t) dt+\frac{\si_1L_0}{8}\int\limits_0^T(\omega_x-l\vph)^2(L_0,t) dt
			\\
			+
			\frac{4k_1}{l^2L_0}\int\limits_0^T(\vph_{x}+\psi+l\omega)^2(0,t) dt+\frac{\si_1L_0}{8}\int\limits_0^T(\omega_x-l\vph)^2(0,t)dt\\
			\frac{\rho_1L_0}{8}\int\limits_0^T  \omega_t^2(L_0,t) dt+\frac{2\rho_1}{l^2L_0}\int\limits_0^T\vph_t^2(L_0,t) dt+
			C(R,T)lot+C(E(0)+E(T)).
		\end{multline}	
		Adding \eqref{ol4} to \eqref{ol3} we obtain
		\begin{multline}
			\label{ol5}
			\frac{\si_1}{4}\int\limits_0^T \int\limits_0^{L_0} (\omega_x-l\vph)^2 dx dt+\frac{\rho_1}{2}\int\limits_0^T \int\limits_0^{L_0} \omega_t^2dx dt +\frac{\rho_1L_0}{8}\int\limits_0^T\omega_t^2(0,t)dt+\frac{\sigma_1L_0}{16} \int\limits_0^T  (\omega_x-l\vph)^2(L_0,t) dt\\+\frac{\sigma_1L_0}{16} \int\limits_0^T  (\omega_x-l\vph)^2(L_0,t) dt\le  \rho_1 \int\limits_0^T \int\limits_0^{L_0} \vph_t^2 dx dt+ k_1(2+17l^2L_0^2)\int\limits_0^T \int\limits_0^{L_0} (\vph_{x}+\psi+l\omega)^2dx dt\\+\frac{4k_1}{l^2L_0}\int\limits_0^T  (\vph_{x}+\psi+l\omega)^2(L_0,t) dt+
			\frac{4k_1}{l^2L_0}\int\limits_0^T(\vph_{x}+\psi+l\omega)^2(0,t) dt\\+\frac{4k_1}{l^2} \int\limits_0^T \int\limits_0^{L_0} \psi_x^2 dx dt
			+\frac{2\rho_1}{l^2L_0}\int\limits_0^T\vph_t^2(L_0,t) dt+C(R,T)lot+C(E(0)+E(T)).
		\end{multline}	
		{\it Step 4.} Now we multiply equation \eqref{eq1} by $-\frac{16}{l^2L_0^2}x\vph_x$ and $-\frac{16}{l^2L_0^2}(x-L_0)\vph_x$ and sum up the results. After integration by parts with respect to $t$ we get
		\begin{multline*}
			\frac{16\rho_1}{l^2L_0^2}\int\limits_0^T \int\limits_0^{L_0}\vph_{t}x\vph_{tx}dx dt+\frac{16\rho_1}{l^2L_0^2}\int\limits_0^T \int\limits_0^{L_0}\vph_{t}(x-L_0)\vph_{tx}dx dt\\
			+\frac{16k_1}{l^2L_0^2}\int\limits_0^T \int\limits_0^{L_0}(\vph_x+\psi+l\om)_xx\vph_xdx dt+\frac{16k_1}{l^2L_0^2}\int\limits_0^T \int\limits_0^{L_0}(\vph_x+\psi+l\om)_x(x-L_0)\vph_xdx dt
		\end{multline*}	
		\begin{multline}
			\label{p1}
			+ \frac{16\si_1}{lL_0^2}\int\limits_0^T \int\limits_0^{L_0}(\om_x-l\vph)x\vph_xdx dt+\frac{16\si_1}{lL_0^2}\int\limits_0^T \int\limits_0^{L_0}(\om_x-l\vph)(x-L_0)\vph_xdx dt \\
			-\frac{16}{l^2L_0^2}\int\limits_0^T \int\limits_0^{L_0}(f_1(\tilde\vph,\tilde\psi,\tilde\omega)-f_1(\hat\vph,\hat\psi,\hat\omega))(2x-L_0)\vph_xdx dt\\
			=\frac{16\rho_1}{l^2L_0^2} \int\limits_0^{L_0}\vph_{t}(x,T)(2x-L_0)\vph_{x}(x,T)dx-\frac{16\rho_1}{l^2L_0^2} \int\limits_0^{L_0}\vph_{t}(x,T)(2x-L_0)\vph_{x}(x,T)dx .
		\end{multline}	
		It is easy to see that
		\begin{equation*}
			\frac{16\rho_1}{l^2L_0^2}\int\limits_0^T \int\limits_0^{L_0}\vph_{t}x\vph_{tx}dx dt+\frac{16\rho_1}{l^2L_0^2}\int\limits_0^T \int\limits_0^{L_0}\vph_{t}(x-L_0)\vph_{tx}dx dt
		\end{equation*}	
		\begin{equation}
			\label{p2}
			\qquad\qquad\qquad\qquad\qquad\qquad\qquad\qquad	=-\frac{16\rho_1}{l^2L_0^2}\int\limits_0^T \int\limits_0^{L_0}\vph_{t}^2 dx dt+\frac{8\rho_1}{l^2L_0}\int\limits_0^T\vph_{t}^2 (L_0,t) dt
		\end{equation}	
		and
		\begin{multline*}
			\frac{16k_1}{l^2L_0^2}\int\limits_0^T \int\limits_0^{L_0}(\vph_x+\psi+l\om)_xx\vph_xdx dt+\frac{16k_1}{l^2L_0^2}\int\limits_0^T \int\limits_0^{L_0}(\vph_x+\psi+l\om)_x(x-L_0)\vph_xdx dt\\
			=-\frac{16k_1}{l^2L_0^2}\int\limits_0^T \int\limits_0^{L_0}(\vph_x+\psi+l\om)^2dx dt+\frac{8k_1}{l^2L_0}\int\limits_0^T (\vph_x+\psi+l\om)^2(0,t) dt\\+\frac{8k_1}{l^2L_0}\int\limits_0^T (\vph_x+\psi+l\om)^2(L_0,t) dt
			-\frac{16k_1}{l^2L_0^2}\int\limits_0^T \int\limits_0^{L_0}(\vph_x+\psi+l\om)_xx(\psi+l\om)dx dt\\
			-\frac{16k_1}{l^2L_0^2}\int\limits_0^T \int\limits_0^{L_0}(\vph_x+\psi+l\om)_x(x-L_0)(\psi+l\om)dx dt\\
			=-\frac{16k_1}{l^2L_0^2}\int\limits_0^T \int\limits_0^{L_0}(\vph_x+\psi+l\om)^2dx dt+\frac{8k_1}{l^2L_0}\int\limits_0^T (\vph_x+\psi+l\om)^2(0,t) dt\\
			+\frac{8k_1}{l^2L_0}\int\limits_0^T (\vph_x+\psi+l\om)^2(L_0,t) dt-\frac{16k_1}{l^2L_0}\int\limits_0^T (\vph_x+\psi+l\om)(L_0,t)(\psi+l\om)(L_0,t) dt
		\end{multline*}	
		\begin{multline}
			\label{p3}
			+\frac{32k_1}{l^2L_0^2}\int\limits_0^T \int\limits_0^{L_0}(\vph_x+\psi+l\om)(\psi+l\om)dx dt+
			+\frac{16k_1}{lL_0^2}\int\limits_0^T \int\limits_0^{L_0}(\vph_x+\psi+l\om)(2x-L_0)(\om_x-l\vph)dx dt\\+\frac{16k_1}{l^2L_0^2}\int\limits_0^T \int\limits_0^{L_0}(\vph_x+\psi+l\om)(2x-L_0)\psi_xdx dt
			+\frac{16k_1}{L_0^2}\int\limits_0^T \int\limits_0^{L_0}(\vph_x+\psi+l\om)(2x-L_0)\vph dx dt.
		\end{multline}	
		Moreover,
		\begin{multline*}
			\frac{16\si_1}{lL_0^2}\int\limits_0^T \int\limits_0^{L_0}(\om_x-l\vph)x\vph_xdx dt+\frac{16\si_1}{lL_0^2}\int\limits_0^T \int\limits_0^{L_0}(\om_x-l\vph)(x-L_0)\vph_xdx dt\\
			=	\frac{16\si_1}{lL_0^2}\int\limits_0^T \int\limits_0^{L_0}(\om_x-l\vph)(2x-L_0)(\vph_x+\psi+l\om)dx dt
		\end{multline*}
		\begin{equation}
			\label{p4}
			\qquad\qquad\qquad	\qquad\qquad\qquad-	\frac{16\si_1}{lL_0^2}\int\limits_0^T \int\limits_0^{L_0}(\om_x-l\vph)(2x-L_0)(\psi+l\om)dx dt.
		\end{equation}	
		Collecting \eqref{p1}--\eqref{p4} and using the estimates
		\begin{multline*}
			\left|\frac{32k_1}{lL_0^2}\int\limits_0^T \int\limits_0^{L_0}(\vph_x+\psi+l\om)(2x-L_0)(\om_x-l\vph)dx dt\right|\\\le \frac{\si_1}{8}\int\limits_0^T \int\limits_0^{L_0}(\om_x-l\vph)^2dx dt+\frac{2046k_1}{l^2L_0^2}\int\limits_0^T \int\limits_0^{L_0}(\vph_x+\psi+l\om)^2dx dt
		\end{multline*}	
		and
		\begin{multline*}
			\left|\frac{16k_1}{l^2L_0^2}\int\limits_0^T \int\limits_0^{L_0}(\vph_x+\psi+l\om)(2x-L_0)\psi_xdx dt\right|\\\le \frac{k_1}{l^2}\int\limits_0^T \int\limits_0^{L_0}\psi_x^2dx dt+\frac{64k_1}{l^2L_0^2}\int\limits_0^T \int\limits_0^{L_0}(\vph_x+\psi+l\om)^2dx dt
		\end{multline*}	
		we come to
		\begin{equation}
			\frac{7k_1}{l^2L_0}\int\limits_0^T (\vph_x+\psi+l\om)^2(L_0,t) dt+\frac{7k_1}{l^2L_0}\int\limits_0^T (\vph_x+\psi+l\om)^2(0,t) dt\qquad\qquad\qquad\qquad\qquad\qquad\qquad\qquad
		\end{equation}	
		\begin{multline}
			\label{ol6}
			+\frac{8\rho_1}{l^2L_0}\int\limits_0^T \vph_t^2(L_0,t) dt\le  \frac{16\rho_1}{l^2L_0^2}\int\limits_0^T \int\limits_0^{L_0} \vph_t^2 dx dt+ \frac{2150 k_1}{l^2L_0^2}\int\limits_0^T \int\limits_0^{L_0} (\vph_{x}+\psi+l\omega)^2dx dt\\+\frac{k_1}{l^2} \int\limits_0^T \int\limits_0^{L_0} \psi_x^2 dx dt+\frac{3\si_1}{16}\int\limits_0^T \int\limits_0^{L_0}(\om_x-l\vph)^2dx dt+C(R,T)lot+C(E(0)+E(T)).
		\end{multline}	
		Adding \eqref{ol6} to \eqref{ol5} we arrive at
		\begin{multline}
			\label{ol7}
			\frac{\si_1}{16}\int\limits_0^T \int\limits_0^{L_0} (\omega_x-l\vph)^2 dx dt+\frac{\rho_1}{2}\int\limits_0^T \int\limits_0^{L_0} \omega_t^2dx dt +\frac{\rho_1L_0}{8}\int\limits_0^T\omega_t^2(L_0,t)dt\\+\frac{\sigma_1L_0}{16} \int\limits_0^T  (\omega_x-l\vph)^2(L_0,t) dt+\frac{\sigma_1L_0}{16} \int\limits_0^T  (\omega_x-l\vph)^2(0,t) dt\\
			+	
			\frac{3k_1}{l^2L_0}\int\limits_0^T (\vph_x+\psi+l\om)^2(L_0,t) dt+\frac{3k_1}{l^2L_0}\int\limits_0^T (\vph_x+\psi+l\om)^2(0,t) dt\\+\frac{6\rho_1}{l^2L_0}\int\limits_0^T \vph_t^2(L_0,t) dt
			\le \rho_1 \left(1+\frac{16}{l^2L_0^2} \right)\int\limits_0^T \int\limits_0^{L_0} \vph_t^2 dx dt\\+ k_1\left(2+17l^2L_0^2+\frac{2150}{l^2L_0^2}\right)\int\limits_0^T \int\limits_0^{L_0} (\vph_{x}+\psi+l\omega)^2dx dt\\+\frac{5k_1}{l^2} \int\limits_0^T \int\limits_0^{L_0} \psi_x^2 dx dt+C(R,T)lot+C(E(0)+E(T)).
		\end{multline}	
		{\it Step 5.} Next we multiply equation \eqref{eq1} by $-\left(1+\frac{18}{l^2L_0^2}\right)\vph$ and integrate by parts with respect to $t$
		\begin{multline}	\label{p5}
			\rho_1\left(1+\frac{18}{l^2L_0^2}\right)\int\limits_0^T \int\limits_0^{L_0}\vph_{t}^2 dx dt+k_1\left(1+\frac{18}{l^2L_0^2}\right)\int\limits_0^T\int\limits_0^{L_0} (\vph_x+\psi+l\om)_x\vph dx dt \\
			+ l\si_1\left(1+\frac{18}{l^2L_0^2}\right)\int\limits_0^T\int\limits_0^{L_0}(\om_x-l\vph)\vph dx dt
			-\left(1+\frac{18}{l^2L_0^2}\right)\int\limits_0^T\int\limits_0^{L_0}(f_1(\tilde\vph, \tilde\psi,\tilde\om)-f_1(\hat\vph,\hat\psi,\hat\om))\vph dx dt=\\
			\rho_1\left(1+\frac{18}{l^2L_0^2}\right) \int\limits_0^{L_0}(\vph_t(x,T)\vph(x,T)-\vph_t(x,0)\vph(x,0))dx.
		\end{multline}	
		Since
		\begin{multline}	\label{p6}
			k_1\left(1+\frac{18}{l^2L_0^2}\right)\int\limits_0^T\int\limits_0^{L_0} (\vph_x+\psi+l\om)_x\vph dx dt
			=-k_1\left(1+\frac{18}{l^2L_0^2}\right)\int\limits_0^T \int\limits_0^{L_0} (\vph_x+\psi+l\om)^2 dx dt\\
			+k_1\left(1+\frac{18}{l^2L_0^2}\right)\int\limits_0^T (\vph_x+\psi+l\om)(L_0,t)\vph(L_0,t)  dt\\
			+k_1\left(1+\frac{18}{l^2L_0^2}\right)\int\limits_0^T (\vph_x+\psi+l\om)(\psi+l\om) dx dt
		\end{multline}	
		we obtain the estimate
		\begin{multline}
			\label{ol8}
			\rho_1\left(1+\frac{17}{l^2L_0^2}\right)\int\limits_0^T \int\limits_0^{L_0}\vph_{t}^2 dx dt\le k_1\left(2+\frac{18}{l^2L_0^2}\right)\int\limits_0^T\int\limits_0^{L_0} (\vph_x+\psi+l\om)^2 dx dt\\
			+\frac{k_1}{l^2L_0}\int\limits_0^T (\vph_x+\psi+l\om)^2(L_0,t) dt+\frac{\si_1}{32}\int\limits_0^T \int\limits_0^{L_0} (\omega_x-l\vph)^2 dx dt\\+C(R,T)lot+C(E(0)+E(T)).	
		\end{multline}	
		Summing up \eqref{ol7} and \eqref{ol8} we get
		\begin{multline*}
			\frac{\si_1}{32}\int\limits_0^T \int\limits_0^{L_0} (\omega_x-l\vph)^2 dx dt+\frac{\rho_1}{2}\int\limits_0^T \int\limits_0^{L_0} \omega_t^2dx dt +\frac{\rho_1L_0}{8}\int\limits_0^T\omega_t^2(L_0,t)dt\\+\frac{\sigma_1L_0}{16} \int\limits_0^T  (\omega_x-l\vph)^2(L_0,t) dt+\frac{\sigma_1L_0}{16} \int\limits_0^T  (\omega_x-l\vph)^2(0,t) dt\\+	
			\frac{2k_1}{l^2L_0}\int\limits_0^T (\vph_x+\psi+l\om)^2(L_0,t) dt+\frac{2k_1}{l^2L_0}\int\limits_0^T (\vph_x+\psi+l\om)^2(0,t) dt\\
			+\frac{6\rho_1}{l^2L_0}\int\limits_0^T \vph_t^2(L_0,t) dt+	\frac{1}{l^2L_0^2}\int\limits_0^T \int\limits_0^{L_0}\vph_{t}^2 dx dt\\
			\le k_1\left(4+17l^2L_0^2+\frac{2200}{l^2L_0^2}\right)\int\limits_0^T (\vph_x+\psi+l\om)^2 dx dt
		\end{multline*}	
		\begin{equation}\label{ol9}
			\qquad\qquad\qquad\qquad\qquad\qquad\qquad\qquad\qquad+\frac{6k_1}{l^2} \int\limits_0^T \int\limits_0^{L_0} \psi_x^2 dx dt+C(R,T)lot+C(E(0)+E(T)).	
		\end{equation}	
		{\it Step 6.} Next we multiply equation \eqref{eq2} by $C_1(\vph_x+\psi+l\om)$ and equation \eqref{eq1} by $C_1\frac{\beta_1}{\rho_1}\psi_x$, where $C_1=2(6+17l^2L_0^2+\frac{2200}{l^2L_0^2})$. Then we sum up the results and integrate by parts with respect to $t$. Taking into account \eqref{c1}, \eqref{c2} we come to
		\begin{multline} \label{p7}
			-\beta_1C_1 \int\limits_0^T \int\limits_0^{L_0}\vph_{t}\psi_{tx}dx dt-\la_1C_1 \int\limits_0^T \int\limits_0^{L_0}(\vph_x+\psi+l\om)_x\psi_x dx dt\\
			- lC_1\la_1\int\limits_0^T \int\limits_0^{L_0}(\om_x-l\vph)\psi_x dx dt +C_1\frac{\beta_1}{\rho_1}\int\limits_0^T \int\limits_0^{L_0}(f_1(\tilde\vph, \tilde\psi, \tilde\om)-f_1(\hat\vph, \hat\psi, \hat\om))\psi_x dx dt\\
			-\be_1C_1 \int\limits_0^T \int\limits_0^{L_0}\psi_{t}(\vph_{xt}+\psi_t+l\om_t)dx dt -\la_1C_1 \int\limits_0^T \int\limits_0^{L_0}\psi_{xx}(\vph_x+\psi+l\om)dx dt\\
			+k_1C_1 \int\limits_0^T \int\limits_0^{L_0}(\vph_x+\psi+l\om)^2 dx dt + C_1 \int\limits_0^T \int\limits_0^{L_0}(\ga(\tilde\psi_t)-\ga(\hat\psi_t)) (\vph_x+\psi+l\om)dx dt\\
			+C_1\int\limits_0^T \int\limits_0^{L_0}(h_1(\tilde\vph, \tilde\psi, \tilde\om)-h_1(\hat\vph, \hat\psi, \hat\om))(\vph_x+\psi+l\om)dx dt
			=\beta_1C_1\int\limits_0^{L_0}\vph_{t}(x,0)\psi_{x}(x,0)dx\\-\beta_1C_1\int\limits_0^{L_0}\vph_{t}(x,T)\psi_{x}(x,T)dx
			+\be_1C_1 \int\limits_0^{L_0}\psi_{t}(x,0)(\vph_{x}+\psi+l\om)(x,0)dx\\
			-\be_1C_1 \int\limits_0^{L_0}\psi_{t}(x,T)(\vph_{x}+\psi+l\om)(x,T)dx.
		\end{multline}	
		Integrating by parts with respect to $x$ we get
		\begin{multline}
			\label{p8}
			\left|\beta_1C_1 \int\limits_0^T \int\limits_0^{L_0}\vph_{t}\psi_{tx}dx dt+\be_1C_1 \int\limits_0^T \int\limits_0^{L_0}\psi_{t}(\vph_{xt}+l\om_t)dx dt\right|\\
			\le \left|\beta_1C_1 \int\limits_0^T \vph_{t}(L_0,t)\psi_{t}(L_0,t) dt+\be_1C_1 l\int\limits_0^T \int\limits_0^{L_0}\psi_{t}\om_tdx dt\right|\le
			\frac{\rho_1}{l^2L_0}\int\limits_0^T \vph_{t}^2(L_0,t) dt\\+	\frac{\beta_1^2C_1^2l^2L_0}{4\rho_1}\int\limits_0^T \psi_{t}^2(L_0,t) dt
			+ \frac{\rho_1}{4}\int\limits_0^T \int\limits_0^{L_0}\om_t^2dx dt+  \frac{\be_1^2C_1^2 l^2}{\rho_1}\int\limits_0^T \int\limits_0^{L_0}\psi_{t}^2dx dt				
		\end{multline}	
		and
		\begin{multline}
			\label{p9}
			\left|	\la_1C_1 \int\limits_0^T \int\limits_0^{L_0}(\vph_x+\psi+l\om)_x\psi_x dx dt+\la_1C_1 \int\limits_0^T \int\limits_0^{L_0}\psi_{xx}(\vph_x+\psi+l\om)dx dt\right|\\=\left|\la_1C_1 \int\limits_0^T (\vph_x+\psi+l\om)(L_0,t)\psi_x(L_0,t) dt-\la_1C_1 \int\limits_0^T (\vph_x+\psi+l\om)(0,t)\psi_x(0,t) dt\right|\\
			\le \frac{k_1}{l^2L_0}\int\limits_0^T (\vph_x+\psi+l\om)^2(L_0,t) dt+\frac{k_1}{l^2L_0}\int\limits_0^T (\vph_x+\psi+l\om)^2(0,t) dt\\
			+\frac{l^2L_0\la_1^2C_1^2}{4k_1}\int\limits_0^T \psi_x^2(L_0,t) dt+\frac{l^2L_0\la_1^2C_1^2}{4k_1}\int\limits_0^T \psi_x^2(0,t) dt.
		\end{multline}	
		Moreover,
		\begin{equation}\label{p10}
			\left|lC_1\la_1\int\limits_0^T \int\limits_0^{L_0}(\om_x-l\vph)\psi_x dx dt \right|\le \frac{\si_1}{64}\int\limits_0^T \int\limits_0^{L_0}(\om_x-l\vph)^2 dx dt+\frac{16l^2C_1^2\la_1^2}{\si_1}\int\limits_0^T \int\limits_0^{L_0}\psi_x^2 dx dt.
		\end{equation}	
		It follows from Lemma \eqref{lem:GammaEst} with $\varepsilon=\frac{k_1 C_1}{4}$
		\begin{multline}\label{p111}
			\left|C_1 \int\limits_0^T \int\limits_0^{L_0} (\ga(\tilde\psi_t)-\ga(\hat\psi_t)) (\vph_x+\psi+l\om)dx dt\right|\\
			\le
			\frac{k_1 C_1}{4} \int\limits_0^T \int\limits_0^{L_0}(\vph_x+\psi+l\om)^2dx dt+C\int\limits_0^T \int\limits_0^{L_0}(\ga(\tilde\psi_t)-\ga(\hat\psi_t))\psi_t dx dt
		\end{multline}	
		Consequently, collecting \eqref{p7}--\eqref{p111} we obtain
		\begin{multline*}
			\frac{C_1k_1}{2}\int\limits_0^T \int\limits_0^{L_0}(\vph_x+\psi+l\om)^2dx dt\le
			\frac{\si_1}{64}\int\limits_0^T \int\limits_0^{L_0}(\om_x-l\vph)^2 dx dt\\
			+\frac{20l^2C_1^2\la_1^2}{\si_1}\int\limits_0^T \int\limits_0^{L_0}\psi_x^2 dx dt+
			C_1\left(\beta_1+\frac{\beta_1^2l^2}{\rho_1}\right)\int\limits_0^T \int\limits_0^{L_0}\psi_t^2 dx dt + \\
			\frac{k_1}{l^2L_0}\int\limits_0^T (\vph_x+\psi+l\om)^2(L_0,t) dt+\frac{k_1}{l^2L_0}\int\limits_0^T (\vph_x+\psi+l\om)^2(0,t) dt
		\end{multline*}	
		\begin{multline}\label{ol10}
			+\frac{l^2L_0\la_1^2C_1^2}{4k_1}\int\limits_0^T \psi_x^2(L_0,t) dt+\frac{l^2L_0\la_1^2C_1^2}{4k_1}\int\limits_0^T \psi_x^2(0,t) dt\\
			+\frac{\rho_1}{l^2L_0}\int\limits_0^T \vph_{t}^2(L_0,t) dt+	\frac{\beta_1^2C_1^2l^2L_0}{4\rho_1}\int\limits_0^T \psi_{t}^2(L_0,t) dt+\frac{\rho_1}{4}\int\limits_0^T \int\limits_0^{L_0}\om_t^2dx dt \\+C\int\limits_0^T \int\limits_0^{L_0}(\ga(\tilde\psi_t)-\ga(\hat\psi_t))\psi_t dx dt+
			C(R,T)lot+C(E(0)+E(T)).
		\end{multline}	
		Combining \eqref{ol10} with \eqref{ol9} we get
		\begin{multline}\label{ol11}
			\frac{\si_1}{64}\int\limits_0^T \int\limits_0^{L_0} (\omega_x-l\vph)^2 dx dt+\frac{\rho_1}{4}\int\limits_0^T \int\limits_0^{L_0} \omega_t^2dx dt +\frac{\rho_1L_0}{8}\int\limits_0^T\omega_t^2(L_0,t)dt\\
			+\frac{\sigma_1L_0}{16} \int\limits_0^T  (\omega_x-l\vph)^2(L_0,t) dt+\frac{\sigma_1L_0}{16} \int\limits_0^T  (\omega_x-l\vph)^2(0,t) dt\\
			+\frac{k_1}{l^2L_0}\int\limits_0^T (\vph_x+\psi+l\om)^2(L_0,t) dt+\frac{k_1}{l^2L_0}\int\limits_0^T (\vph_x+\psi+l\om)^2(0,t) dt\\
			+\frac{5\rho_1}{l^2L_0}\int\limits_0^T \vph_t^2(L_0,t) dt+	\frac{1}{l^2L_0^2}\int\limits_0^T \int\limits_0^{L_0}\vph_{t}^2 dx dt+ 2k_1\int\limits_0^T \int\limits_0^{L_0} (\vph_x+\psi+l\om)^2 dx dt\\
			\le \left(\frac{6k_1}{l^2}+\frac{20l^2C_1^2\lambda_1^2}{\sigma_1}\right) \int\limits_0^T \int\limits_0^{L_0} \psi_x^2 dx dt+C_1\left(\beta_1+\frac{\beta_1^2l^2}{\rho_1}\right)\int\limits_0^T \int\limits_0^{L_0}\psi_t^2 dx dt\\
			+\frac{l^2L_0\la_1^2C_1^2}{4k_1}\int\limits_0^T \psi_x^2(L_0,t) dt+\frac{l^2L_0\la_1^2C_1^2}{4k_1}\int\limits_0^T \psi_x^2(0,t) dt
			+\frac{\beta_1^2C_1^2l^2L_0}{4}\int\limits_0^T \psi_{t}^2(L_0,t) dt\\+C\int\limits_0^T \int\limits_0^{L_0}(\ga(\tilde\psi_t)-\ga(\hat\psi_t))\psi_t dx dt+C(R,T)lot+C(E(0)+E(T)).	
		\end{multline}	
		{\it Step 7.} Our next step is to multiply equation \eqref{eq2} by $-C_2x\psi_x-C_2(x-L_0)\psi_x$, where $C_2=\frac{l^2\la_1C_1^2}{k_1}$. After integration by parts with respect to $t$ we obtain
		\begin{equation*} 	
			\be_1C_2\int\limits_0^T \int\limits_0^{L_0}\psi_{t}x \psi_{xt} dx dt+\be_1C_2\int\limits_0^T \int\limits_0^{L_0}\psi_{t}(x-L_0) \psi_{xt} dx dt\qquad\qquad\qquad\qquad\qquad\qquad
		\end{equation*}	
		\begin{multline} 	\label{p11}
			+\la_1C_2 \int\limits_0^T \int\limits_0^{L_0}\psi_{xx}x \psi_{x} dx dt+\la_1C_2 \int\limits_0^T \int\limits_0^{L_0}\psi_{xx}(x-L_0) \psi_{x} dx dt\\
			-k_1C_2\int\limits_0^T \int\limits_0^{L_0}(\vph_x+\psi+l\om)(2x-L_0)\psi_{x} dx dt
			- C_2 \int\limits_0^T \int\limits_0^{L_0}(\gamma(\tilde\psi_t)-\gamma(\hat\psi_t))(2x-L_0)\psi_{x} dx dt\\
			+ \int\limits_0^T \int\limits_0^{L_0}(h_1(\tilde\vph,\tilde\psi, \tilde\om)-h_1(\hat\vph,\hat\psi, \hat\om))(2x-L_0)\psi_{x} dx dt\\
			= \be_1C_2 \int\limits_0^{L_0}\psi_{t}(x,T)(2x-L_0) \psi_{x}(x,T) dx -\be_1C_2 \int\limits_0^{L_0}\psi_{t}(x,0)(2x-L_0) \psi_{x}(x,0) dx.
		\end{multline}	
		After integration by parts with respect to $x$ we get
		\begin{multline}
			\label{p12}
			\be_1C_2\int\limits_0^T \int\limits_0^{L_0}\psi_{t}x \psi_{xt} dx dt+\be_1C_2\int\limits_0^T \int\limits_0^{L_0}\psi_{t}(x-L_0) \psi_{xt} dx dt\\=-\be_1C_2\int\limits_0^T \int\limits_0^{L_0}\psi_{t}^2 dx dt+\frac{\be_1C_2L_0}{2}\int\limits_0^T \psi_{t}^2(L_0,t)dt
		\end{multline}	
		and
		\begin{multline}
			\label{p13}
			\la_1C_2 \int\limits_0^T \int\limits_0^{L_0}\psi_{xx}x \psi_{x} dx dt+\la_1C_2 \int\limits_0^T \int\limits_0^{L_0}\psi_{xx}(x-L_0) \psi_{x} dx dt\\
			=\frac{\la_1C_2L_0}{2} \int\limits_0^T\psi_{x}^2(L_0,t) dt+\frac{\la_1C_2L_0}{2} \int\limits_0^T\psi_{x}^2(0,t) dt-
			\la_1C_2 \int\limits_0^T \int\limits_0^{L_0}\psi_{x}^2 dx dt.
		\end{multline}	
		Furthermore,
		\begin{multline}
			\label{p14}
			\left| k_1C_2\int\limits_0^T \int\limits_0^{L_0}(\vph_x+\psi+l\om)(2x-L_0)\psi_{x} dx dt\right|\\
			\le
			k_1\int\limits_0^T \int\limits_0^{L_0}(\vph_x+\psi+l\om)^2 dx dt+
			\frac{k_1C_2^2L_0^2}{4}\int\limits_0^T \int\limits_0^{L_0}\psi_{x}^2 dx dt.
		\end{multline}	
		By Lemma \eqref{lem:GammaEst} with $\varepsilon=\frac{k_1 C2^2L_0^2}{4}$ we have
		\begin{equation}
			\label{p141}
			\left| C_2 \int\limits_0^T \int\limits_0^{L_0}\psi_t(2x-L_0)\psi_{x} dx dt\right|\le
			\frac{k_1C_2^2L_0^2}{4}\int\limits_0^T \int\limits_0^{L_0}\psi_{x}^2 dx dt+C \int\limits_0^T \int\limits_0^{L_0}(\gamma(\tilde\psi_t)-\gamma(\hat\psi_t))\psi_tdx dt.
		\end{equation}	
		As a result of \eqref{p11}-- \eqref{p141} we obtain the estimate
		\begin{multline}
			\label{ol12}
			\frac{\be_1C_2L_0}{2}\int\limits_0^T \psi_{t}^2(L_0,t)dt+\frac{\la_1C_2L_0}{2} \int\limits_0^T\psi_{x}^2(L_0,t) dt+\frac{\la_1C_2L_0}{2} \int\limits_0^T\psi_{x}^2(0,t)dt\\
			\le k_1 \int\limits_0^T \int\limits_0^{L_0}(\vph_x+\psi+l\om)^2 dx dt+
			\left(k_1C_2^2L_0^2+\lambda_1C_2\right)\int\limits_0^T \int\limits_0^{L_0}\psi_{x}^2 dx dt+\beta_1C_2 \int\limits_0^T \int\limits_0^{L_0}\psi_t^2dx dt\\
			+C \int\limits_0^T \int\limits_0^{L_0}(\gamma(\tilde\psi_t)-\gamma(\hat\psi_t))\psi_tdx dt+C(R,T)lot+C(E(0)+E(T)).
		\end{multline}
		Summing up \eqref{ol11} and \eqref{ol12} and using \eqref{c2} we infer
		\begin{multline*}
			\frac{\si_1}{64}\int\limits_0^T \int\limits_0^{L_0} (\omega_x-l\vph)^2 dx dt+\frac{\rho_1}{4}\int\limits_0^T \int\limits_0^{L_0} \omega_t^2dx dt +\frac{\rho_1L_0}{8}\int\limits_0^T\omega_t^2(L_0,t)dt\\
			+\frac{\sigma_1L_0}{16} \int\limits_0^T  (\omega_x-l\vph)^2(L_0,t) dt+\frac{\sigma_1L_0}{16} \int\limits_0^T  (\omega_x-l\vph)^2(0,t) dt\\
			+\frac{k_1}{l^2L_0}\int\limits_0^T (\vph_x+\psi+l\om)^2(L_0,t) dt+\frac{k_1}{l^2L_0}\int\limits_0^T (\vph_x+\psi+l\om)^2(0,t) dt\\
			+\frac{5\rho_1}{l^2L_0}\int\limits_0^T \vph_t^2(L_0,t) dt+	\frac{1}{l^2L_0^2}\int\limits_0^T \int\limits_0^{L_0}\vph_{t}^2 dx dt+ k_1\int\limits_0^T \int\limits_0^{L_0}(\vph_x+\psi+l\om)^2 dx dt\\
			\frac{l^2L_0\la_1^2C_1^2}{4k_1}\int\limits_0^T \psi_x^2(L_0,t) dt+\frac{l^2L_0\la_1^2C_1^2}{4k_1}\int\limits_0^T \psi_x^2(0,t) dt\\
			+\frac{\beta_1^2C_1^2l^2L_0}{4\rho_1}\int\limits_0^T \psi_{t}^2(L_0,t) dt\le \left(\frac{6k_1}{l^2}+\frac{20l^2C_1^2\lambda_1^2}{\sigma_1}+\lambda_1 C_2+k_1C_2^2L_0^2\right) \int\limits_0^T \int\limits_0^{L_0} \psi_x^2 dx dt
		\end{multline*}	
		\begin{multline} \label{ol13}
			+\left((C_1+C_2)\beta_1+\frac{C_1\beta_1^2l^2}{\rho_1}\right)\int\limits_0^T \int\limits_0^{L_0}\psi_t^2 dx dt\\+C \int\limits_0^T \int\limits_0^{L_0}(\gamma(\tilde\psi_t)-\gamma(\hat\psi_t))\psi_tdx dt
			+C(R,T)lot+C(E(0)+E(T)).		
		\end{multline}	
		{\it Step 8.}	Now we multiply equation \eqref{eq2} by $C_3\psi$, where $C_3=\frac{2}{\lambda_1}\left(\frac{6k_1}{l^2}+\frac{20l^2C_1^2\lambda_1^2}{\sigma_1}+\lambda_1 C_2+k_1C_2^2L_0^2\right)$ and integrate by parts with respect to $t$
		\begin{multline}
			-C_3\be_1\int\limits_0^T \int\limits_0^{L_0}\psi_{t}^2 dx dt-\la_1C_3 \int\limits_0^T \int\limits_0^{L_0}\psi_{xx}\psi dx dt +k_1C_3\int\limits_0^T \int\limits_0^{L_0}(\vph_x+\psi+l\om)\psi dx dt\\
			+ C_3\int\limits_0^T \int\limits_0^{L_0} (\ga(\tilde\psi_t)-\ga(\hat\psi_t))\psi dx dt  +C_3\int\limits_0^T \int\limits_0^{L_0}(h_1(\tilde\vph,\tilde\psi,\tilde\om)-h_1(\hat\vph,\hat\psi,\hat\om))\psi dx dt
		\end{multline}	
		\begin{equation}
			\label{p15}
			\qquad\qquad\qquad\qquad=
			C_3\be_1\int\limits_0^{L_0}\psi_{t}(x,0)\psi(x,0) dx-
			C_3\be_1\int\limits_0^{L_0}\psi_{t}(x,T)\psi(x,T) dx
		\end{equation}	
		After integration by parts we infer the estimate
		\begin{multline}
			\label{ol14}
			\la_1C_3 \int\limits_0^T \int\limits_0^{L_0}\psi_x^2 dx dt \le
			\frac{k_1}{2}\int\limits_0^T (\vph_x+\psi+l\om)^2 dx dt+C_3\beta_1\int\limits_0^T \int\limits_0^{L_0} \psi_t^2 dx dt
			+	\frac{l^2L_0\la_1^2C_1^2}{8k_1}\int\limits_0^T \psi_x^2(L_0,t) dt
			\\
			+C \int\limits_0^T \int\limits_0^{L_0}(\gamma(\tilde\psi_t)-\gamma(\hat\psi_t))\psi_tdx dt	+C(R,T)lot+C(E(0)+E(T)).
		\end{multline}
		Combining \eqref{ol14} with \eqref{ol13} we obtain
		\begin{multline*}
			\frac{\si_1}{64}\int\limits_0^T \int\limits_0^{L_0} (\omega_x-l\vph)^2 dx dt+\frac{\rho_1}{4}\int\limits_0^T \int\limits_0^{L_0} \omega_t^2dx dt +\frac{\rho_1L_0}{8}\int\limits_0^T\omega_t^2(L_0,t)dt\\
			+\frac{\sigma_1L_0}{16} \int\limits_0^T  (\omega_x-l\vph)^2(L_0,t) dt+\frac{\sigma_1L_0}{16} \int\limits_0^T  (\omega_x-l\vph)^2(0,t) dt\\
			+\frac{k_1}{l^2L_0}\int\limits_0^T (\vph_x+\psi+l\om)^2(L_0,t) dt+\frac{k_1}{l^2L_0}\int\limits_0^T (\vph_x+\psi+l\om)^2(0,t) dt
		\end{multline*}
		\begin{multline}\label{ol15}
			+\frac{5\rho_1}{l^2L_0}\int\limits_0^T \vph_t^2(L_0,t) dt+	\frac{1}{l^2L_0^2}\int\limits_0^T \int\limits_0^{L_0}\vph_{t}^2 dx dt+ \frac{k_1}{2}\int\limits_0^T \int\limits_0^{L_0} (\vph_x+\psi+l\om)^2 dx dt\\
			\frac{l^2L_0\la_1^2C_1^2}{8k_1}\int\limits_0^T \psi_x^2(L_0,t) dt+\frac{l^2L_0\la_1^2C_1^2}{4k_1}\int\limits_0^T \psi_x^2(0,t) dt
			+\frac{\beta_1^2C_1^2l^2L_0}{4\rho_1}\int\limits_0^T \psi_{t}^2(L_0,t) dt\\
			+ \left(\frac{6k_1}{l^2}+\frac{20l^2C_1^2\lambda_1^2}{\sigma_1}+\lambda_1 C_2+k_1C_2^2L_0^2\right) \int\limits_0^T \int\limits_0^{L_0} \psi_x^2 dx dt\\
			\le\left((C_1+C_2)\beta_1+\frac{C_1\beta_1^2l^2}{\rho_1}+C_3\beta_1\right)\int\limits_0^T \int\limits_0^{L_0}\psi_t^2 dx dt\\
			+C \int\limits_0^T \int\limits_0^{L_0}(\gamma(\tilde\psi_t)-\gamma(\hat\psi_t))\psi_tdx dt	+C(R,T)lot+C(E(0)+E(T)).	
		\end{multline}	
		{\it Step 9.} Consequently, it follows from \eqref{ol15} and assumption \eqref{GammaLip2} for any $l>0$ where exist  constants $M_i$, $i=\overline{\{1,3\}}$ (depending on $l$) such that
		\begin{multline}
			\label{l1}
			\int\limits_0^T E_1(t) dt+\int\limits_0^T B_1(t) dt\le M_1 \int\limits_0^T \int\limits_0^{L_0}(\gamma(\tilde\psi_t)-\gamma(\hat\psi_t))\psi_tdx dt\\+M_2(R,T)lot+M_3(E(T)+E(0)),
		\end{multline}	
		where
		\begin{multline}
			\label{b1}
			B_1(t)=\int\limits_0^T  (\omega_x-l\vph)^2(L_0,t) dt+	
			\int\limits_0^T (\vph_x+\psi+l\om)^2(L_0,t) dt+\int\limits_0^T \psi_x^2(L_0,t) dt \\+\int\limits_0^T\omega_t^2(L_0,t)dt
			+\int\limits_0^T \psi_{t}^2(L_0,t) dt+\int\limits_0^T \vph_t^2(L_0,t) dt.
		\end{multline}	
		{\it Step 10.} Finally, we multiply equation \eqref{eq4} by $(x-L)u_x$, equation \eqref{eq5} by $(x-L)v_x$, and \eqref{eq6} by $(x-L)w_x$. Summing up the results and integrating by parts with respect to $t$ we arrive at
		\begin{multline*}
			- \rho_2\int\limits_0^T \int\limits_{L_0}^Lu_{t}(x-L)u_{tx} dx dt-k_2\int\limits_0^T \int\limits_{L_0}^L(u_x+v+lw)_x (x-L)u_{x} dx dt\\
			- l\si_2\int\limits_0^T \int\limits_{L_0}^L(w_x-lu)(x-L)u_{x} dx dt +\int\limits_0^T \int\limits_{L_0}^L(f_2(\tilde u, \tilde v, \tilde w)-f_2(\hat u, \hat v, \hat w))(x-L)u_{x} dx dt
		\end{multline*}	
		\begin{multline}\label{r1}
			-\be_2\int\limits_0^T \int\limits_{L_0}^L v_{t} (x-L)v_{xt} dx dt-\la_2 \int\limits_0^T \int\limits_{L_0}^L v_{xx}(x-L)v_{x} dx dt\\
			+k_2 \int\limits_0^T \int\limits_{L_0}^L(u_x+v+lw)(x-L)v_{x} dx dt +\int\limits_0^T \int\limits_{L_0}^L(h_2(\tilde u, \tilde v, \tilde w)-h_2(\hat u, \hat v, \hat w))(x-L)v_{x} dx dt\\
			-\rho_2\int\limits_0^T \int\limits_{L_0}^Lw_{t}(x-L)w_{xt} dx dt- \si_2\int\limits_0^T \int\limits_{L_0}^L(w_x-lu)_x(x-L)w_{x} dx dt\\
			+lk_2\int\limits_0^T \int\limits_{L_0}^L(u_x+v+lw)(x-L)w_{x} dx dt+\int\limits_0^T \int\limits_{L_0}^L(g_2(\tilde u, \tilde v, \tilde w)-g_2(\hat u, \hat v, \hat w))(x-L)w_{x} dx dt=\\
			- \rho_2 \int\limits_{L_0}^L (x-L)((u_{t}u_{x})(x,T) - (u_{t}u_{x})(x,0))dx
			-\be_2 \int\limits_{L_0}^L (x-L)((v_{t}v_{x})(x,T) -(v_{t}v_{x})(x,0))dx \\
			- \rho_2 \int\limits_{L_0}^L (x-L)((w_{t}w_{x})(x,T) - (w_{t}w_{x})(x,0))dx .
		\end{multline}	
		After integration by parts with respect to $x$ we infer
		\begin{multline}
			\label{r2}
			- \rho_2\int\limits_0^T \int\limits_{L_0}^Lu_{t}(x-L)u_{tx} dx-
			\be_2\int\limits_0^T \int\limits_{L_0}^L v_{t} (x-L)v_{xt} dx dt
			-\rho_2\int\limits_0^T \int\limits_{L_0}^Lw_{t}(x-L)w_{xt} dx dt\\=
			\frac{\rho_2}{2}\int\limits_0^T \int\limits_{L_0}^L u_{t}^2 dx+
			\frac{\be_2}{2}\int\limits_0^T \int\limits_{L_0}^L v_{t}^2 dx dt
			+\frac{\rho_2}{2}\int\limits_0^T \int\limits_{L_0}^Lw_{t}^2 dx dt\\-
			\frac{\rho_2(L-L_0)}{2}\int\limits_0^T u_{t}^2(L_0) dt-
			\frac{\be_2(L-L_0)}{2}\int\limits_0^T  v_{t}^2(L_0)  dt
			-\frac{\rho_2(L-L_0)}{2}\int\limits_0^T w_{t}^2(L_0) dt
		\end{multline}	
		and
		\begin{multline*}
			-k_2\int\limits_0^T \int\limits_{L_0}^L(u_x+v+lw)_x (x-L)u_{x} dx dt- l\si_2\int\limits_0^T \int\limits_{L_0}^L(w_x-lu)(x-L)u_{x} dx dt\\
			-\la_2 \int\limits_0^T \int\limits_{L_0}^L v_{xx}(x-L)v_{x} dx dt+k_2 \int\limits_0^T \int\limits_{L_0}^L(u_x+v+lw)(x-L)v_{x} dx dt
		\end{multline*}	
		\begin{multline}
			\label{r3}
			- \si_2\int\limits_0^T \int\limits_{L_0}^L(w_x-lu)_x(x-L)w_{x} dx dt+lk_2\int\limits_0^T \int\limits_{L_0}^L(u_x+v+lw)(x-L)w_{x} dx dt=\\
			-k_2\int\limits_0^T \int\limits_{L_0}^L(u_x+v+lw)_x (x-L)(u_x+v+lw) dx dt\\-\si_2  \int\limits_0^T \int\limits_{L_0}^L(w_x-lu)_x(x-L)(w_x-lu) dx dt -\la_2 \int\limits_0^T \int\limits_{L_0}^L v_{xx}(x-L)v_{x} dx dt\\
			- l\si_2(L-L_0)\int\limits_0^T(w_x-lu)(L_0)u(L_0)dt
			+k_2(L-L_0) \int\limits_0^T(u_x+v+lw)(L_0)v(L_0)dt\\
			+lk_2(L-L_0)\int\limits_0^T (u_x+v+lw)(L_0)w(L_0)dt=\\
			-\frac{k_2(L-L_0)}{2}\int\limits_0^T (u_x+v+lw)^2(L_0) dt+\frac{k_2}{2}\int\limits_0^T \int\limits_{L_0}^L(u_x+v+lw)^2 dx dt\\
			+\frac{\si_2}{2}  \int\limits_0^T \int\limits_{L_0}^L(w_x-lu)^2 dx dt -\frac{\si_2(L-L_0)}{2}  \int\limits_0^T(w_x-lu)^2(L_0) dt\\
			+\frac{\la_2}{2} \int\limits_0^T \int\limits_{L_0}^L v_{x}^2 dx dt-\frac{\la_2(L-L_0) }{2}\int\limits_0^T v_{x}^2(L_0) dt
			- l\si_2(L-L_0)\int\limits_0^T(w_x-lu)(L_0)u(L_0)dt\\
			+k_2(L-L_0) \int\limits_0^T(u_x+v+lw)(L_0)v(L_0)dt\\
			+lk_2(L-L_0)\int\limits_0^T (u_x+v+lw)(L_0)w(L_0)dt.
		\end{multline}	
		Consequently, it follows from \eqref{r1}--\eqref{r3} that  for any $l>0$ where exist  constants $M_4, M_5, M_6>0$  such that
		\begin{equation}
			\label{l2}
			\int\limits_0^T E_2(t) dt\le M_4\int\limits_0^T B_2(t) dt+ M_5(R,T)lot+M_6(E(T)+E(0)),
		\end{equation}	
		where
		\begin{multline}
			\label{b2}
			B_2(t)=\int\limits_0^T  (w_x-lu)^2(L_0,t) dt+	
			\int\limits_0^T (u_x+v+lw)^2(L_0,t) dt+\int\limits_0^T v_x^2(L_0,t) dt \\+\int\limits_0^Tw_t^2(L_0,t)dt
			+\int\limits_0^T v_{t}^2(L_0,t) dt+\int\limits_0^T u_t^2(L_0,t) dt.
		\end{multline}	
		Then, due to transmission conditions \eqref{TC1}--\eqref{TC4}  there exist $\delta, M_7, M_8>0$ (depending on $l$), such that
		\begin{equation}
			\label{l3}
			\int\limits_0^T E(t) dt\le \delta \int\limits_0^T \int\limits_0^{L_0}(\gamma(\tilde\psi_t)-\gamma(\hat\psi_t))\psi_t dx dt+  M_7(R,T)lot+M_8(E(T)+E(0)).
		\end{equation}	
		It follows from \eqref{En1}  that there exists $C>0$ such that
		\begin{equation}
			\label{En11}
			\int\limits_0^T \int\limits_0^{L_0}(\gamma(\tilde\psi_t)-\gamma(\hat\psi_t)) \psi_t dx dt\le C\left( E(0)+ \int\limits_0^T | H(\hat U(t),\tilde U(t))|  dt\right).
		\end{equation}
		By Lemma \ref{lem:lip} we have that for any $\varepsilon>0$ there exists $C(\varepsilon, R)>0$ such that
		\begin{equation} 	\label{HEst}
			\int\limits_0^T | H(\hat U(t), \tilde U(t))|  dt\le \varepsilon\int\limits_0^T\int\limits_0^{L_0} E(t)  dx dt+ C(\varepsilon, R,T)lot.
		\end{equation}
		Combining \eqref{HEst} with \eqref{En11} we arrive at
		\begin{equation}
			\label{DampEst}
			\int\limits_0^T \int\limits_0^{L_0}(\gamma(\tilde\psi_t)-\gamma(\hat\psi_t)) \psi_t dx dt\le C E(0)+C(R,T)lot .
		\end{equation}
		Substituting \eqref{DampEst} into \eqref{l3} we obtain
		\begin{equation}
			\label{l4}
			\int\limits_0^T E(t) dt\le  C(R,T)lot+C(E(T)+E(0))
		\end{equation}	
		for some $C, C(R,T)>0$.\\
		Our remaining task is to estimate the last term in \eqref{En2}.
		\begin{equation}\label{En22}
			\left|	\int\limits_0^T \int\limits_t^T  H(\hat U(s),\tilde U(s)) ds dt\right|\le  \int\limits_0^T E(t) dt+T^3C(R) lot.
		\end{equation}
		Then, it follows from \eqref{En2} and \eqref{En22} that
		\begin{equation}\label{En222}
			TE(T)\le  C\int\limits_0^T E(t) dt+C(T,R) lot.
		\end{equation}
		Then the combination of \eqref{En222} with \eqref{l4} leads to
		\begin{equation}
			\label{l}
			TE(T)\le  C(R,T)lot+C(E(T)+E(0)).
		\end{equation}
		Choosing $T$ large enough one can obtain estimate \eqref{te} which together with Theorem \ref{theoremCL} immediately leads to the asymptotic smoothness of the system.
	\end{proof}
	\subsection{Existence of attractors.}
	The following statement collects criteria on existence and properties of attractors to gradient systems.
	\begin{theorem}[\cite{CFR, CL}]
		\label{abs}
		Assume that $( H, S_t)$ is a gradient asymptotically smooth  dynamical
		system. Assume its Lyapunov function $L(y)$ is bounded from above on any
		bounded subset of $H$ and the set $\cW_R=\{y: L(y) \le R\}$ is bounded for every $R$. If the
		set $\EuScript N$ of stationary points of $(H, S_t)$ is bounded, then $(S_t, H)$ possesses a compact
		global attractor. Moreover,
		the global attractor  consists of  full trajectories
		$\gamma=\{ U(t)\, :\, t\in\R\}$ such that
		\begin{equation}\label{conv-N}
			\lim_{t\to -\infty}{\rm dist}_{H}(U(t),\EuScript N)=0 ~~
			\mbox{and} ~~ \lim_{t\to +\infty}{\rm dist}_{H}(U(t),\EuScript N)=0
		\end{equation}
		and
		\begin{equation}\label{7.4.1}
			\lim_{t\to +\infty}{\rm dist}_{H}(S_tx,\EuScript N)=0
			~~\mbox{for any $x \in H$;}
		\end{equation}
		that is, any trajectory stabilizes to the set $\EuScript N$ of  stationary points.
	\end{theorem}
	Now we state the result on the existence of an attractor.
	\begin{theorem}
		\label{th:attr}
		Let assumptions of Theorems \ref{th:grad}, \ref{th:AsSmooth},  hold true, moreover,
		\begin{gather*}
			\liminf\limits_{|s|\to \infty}\frac{h_1(s)}{s}>0, \tag{N5}\label{fl}\\
			\nabla\cF_2(u,v,w)(u,v,w)-a_1\cF_2(u,v,w)\ge -a_2,\qquad  a_i\ge 0.
		\end{gather*}
		Then, the dynamical system $(H, S_t)$ generated by \eqref{Eq1}-\eqref{TC4} possesses a compact global attractor $\mathfrak A$ possessing properties \eqref{conv-N}, \eqref{7.4.1}.
	\end{theorem}
	\begin{proof}
		In view  of Theorems \ref{th:grad}, \ref{th:AsSmooth}, \ref{abs} our remaining task is to show the boundedness of the set of stationary points and the set $W_R=\{Z:  L(Z)\le R\}$, where $L$ is given by \eqref{lap}.\par
		The second statement follows immediately from the structure of function $ L$ and property \eqref{fl}.\par
		The first statement can be easily shown by energy-like estimates  for stationary solutions taking into account \eqref{fl}.
	\end{proof}
	\section{Singular Limits on finite time intervals}
	\subsection{Singular limit $l\arr 0$}
	Let the nonlinearities $f_j,h_j, g_j$ are such that
	\begin{align*}
		& f_1(\vph,\psi,\om)=f_1(\vph,\psi), &&   h_1(\vph,\psi,\om)=h_1(\vph,\psi), &&  g_1(\vph,\psi,\om)=g_1(\om), \\
		& f_2(u,v,w)=f_2(u,v), &&   h_2(u,v,w)=h_2(u,v), &&  g_2(u,v,w)=g_2(w). \tag{N6}\label{Ndecouple}
	\end{align*}
	If we formally set $l=0$ in  \eqref{AEq}-\eqref{AIC}, we obtain the contact problem for a straight Timoshenko beam
	\begin{align}
		& \rho_1\vph_{tt}-k_1(\vph_x+\psi)_x +{f_1(\vph,\psi)}=p_1(x,{t}), &(x,t)\in (0,L_0)\times (0,T), \label{TimEq1}\\
		& \beta_1\psi_{tt} -\la_1 \psi_{xx} +k_1(\vph_x+\psi) +\ga(\psi_t) +h_1(\vph,\psi)=r_1(x,{t}), &(x,t)\in (0,L_0)\times (0,T),   \\
		& \rho_2u_{tt}-k_2(u_x+v)_x  +f_2(u,v)=p_2(x,{t}), &(x,t)\in (L_0,L)\times (0,T),  \\
		& \beta_2v_{tt} -\la_2 v_{xx} +k_2(u_x+v) +h_2(u,v)=r_2(x,{t}), &(x,t)\in (L_0,L)\times (0,T), \\
		&\vph(0,t)=\psi(0,t)=0, \quad u(L,t)=v(L,t)=0,\\
		& \vph(L_0,t)=u(L_0,t), \quad \psi(L_0,t)=v(L_0,t),\\
		& k_1(\vph_x+\psi)(L_0,t)=k_2(u_x+v)(L_0,t), \, \la_1 \psi_{x}(L_0,t)=\la_2 v_{x}(L_0,t),\label{TimEq2}
	\end{align}
	and an independent contact problem for wave equations
	\begin{align}
		& \rho_1\om_{tt}- \si_1\om_{xx}+{g_1(\om)}=q_1(x,t), &(x,t)\in (0,L_0)\times (0,T), \label{WaveEq1}\\
		& \rho_2w_{tt}- \si_2 w_{xx}+g_2(w)=q_2(x,{t}), &(x,t)\in (L_0,L)\times (0,T), \\
		& \si_1 \om_x(L_0,t)=\si_2w_x(L_0,t),\quad \om(L_0,t)=w(L_0,t), \\
		&w(L,t)=0, \quad \om(0,t)=0. \label{WaveEq2}
	\end{align}
	The following theorem gives an answer, how close are solutions to \eqref{AEq}-\eqref{AIC}  to the solution of decoupled system \eqref{TimEq1}-\eqref{WaveEq2} when $l\arr 0$.
	\begin{theorem}
		Assume that the conditions of Theorem \ref{th:WeakWP}, \eqref{GammaLip1} and \eqref{Ndecouple} hold.
		Let $\Phi^{(l)}$ be the solution to \eqref{AEq}-\eqref{AIC} with the fixed $l$ and the initial data
		\begin{equation*}
			\Phi(x,0)=(\vph_0,\psi_0,\om_0, u_0,v_0,w_0)(x), \quad 	\Phi_t(x,0)=(\vph_1,\psi_1,\om_1, u_1,v_1,w_1)(x).
		\end{equation*}
		Then for every $T>0$
		\begin{align*}
			&\Phi^{(l)}  \stackrel{\ast}{\rightharpoonup} (\vph,\psi,\om, u,v,w) \quad &\mbox{in } L^\infty(0,T;H_d) \; &\mbox{ as } l\arr 0,\\
			&\Phi^{(l)}_t  \stackrel{\ast}{\rightharpoonup} (\vph_t,\psi_t,\om_t, u_t,v_t,w_t) \quad &\mbox{in } L^\infty(0,T;H_v)\; &\mbox{ as } l\arr 0,
		\end{align*}
		where $(\vph,\psi, u,v)$ is the solution to \eqref{TimEq1}-\eqref{TimEq2} with the initial conditions
		\begin{equation*}
			(\vph,\psi, u,v)(x,0)=(\vph_0,\psi_0, u_0,v_0)(x), \quad 	(\vph_t,\psi_t, u_t,v_t)(x,0)=(\vph_1,\psi_1, u_1,v_1)(x),
		\end{equation*}
		and $(\om, w)$ is the solution to \eqref{WaveEq1}-\eqref{WaveEq2} with the initial conditions
		\begin{equation*}
			(\om,w)(x,0)=(\om_0, w_0)(x), \quad 	(\om_t,w_t)(x,0)=(\om_1, w_1)(x).
		\end{equation*}
	\end{theorem}
	
	The proof is similar to that of Theorem 3.1 \cite{MaMo2017} for the homogeneous Bresse beam with obvious changes, except for the limit transition in the nonlinear dissipation term. For the future use we formulate it as a lemma.
	\begin{lemma}\label{lem:DissTrans}
		Let \eqref{GammaLip1} holds. Then
		\begin{equation*}
			\int_0^T \int_0^{L_0} \ga(\psi^{(l)}(x,t))\ga^1(x,t) dxdt \arr \int_0^T \int_0^{L_0} \ga(\psi(x,t))\ga^1(x,t) dxdt  \quad \mbox{ as } l\arr 0
		\end{equation*}
		for every $\ga^1\in L^2(0,T;H^1(0,L_0))$.
	\end{lemma}
	\begin{proof}
		Since \eqref{DCont} and  \eqref{GammaLip1}   hold $|\ga(s)|\le Ms$, therefore 
		\begin{equation*}
			||\ga(\psi^{(l)})||_{L^\infty(0,T; L^2(0,L_0))} \le C( ||\psi^{(l)}||_{L^\infty(0,T; L^2(0,L_0))}).
		\end{equation*}
		Thus, due to Lemmas \ref{lem:ASelfAdjont}, \ref{lem:lip} the sequence
		\begin{equation*}
			R\Phi^{(l)}_{tt}=A\Phi^{(l)}+\Ga(\Phi^{(l)}_t) + F(\Phi^{(l)}) + P
		\end{equation*}
		is bounded in $L^\infty(0,T; H^{-1}(0,L))$ and we can extract from $\Phi^{(l)}_{tt}$ a subsequence, that converges $\ast$-weakly in  $L^\infty(0,T; H^{-1}(0,L))$. Thus, 
		\begin{equation*}
			\Phi^{(l)}_t \arr \Phi_t \quad \mbox{strongly in } L^2(0,T; H^{-\ep}(0,L)), \; \ep>0.
		\end{equation*}
		Consequently,
		\begin{multline*}
			\left|\int_0^T \int_0^{L_0} (\ga(\psi^{(l)}(x,t))-\ga(\psi(x,t))) \ga^1(x,t) dxdt \right| \le \\
			C(L) \int_0^T \int_0^{L_0}  |\psi^{(l)}(x,t)-\psi(x,t)| |\ga^1(x,t)|dxdt \arr 0.
		\end{multline*}
	\end{proof}

\begin{figure}[h!]
		\centering
		\includegraphics[width=0.9\textwidth, height=0.27\textheight]{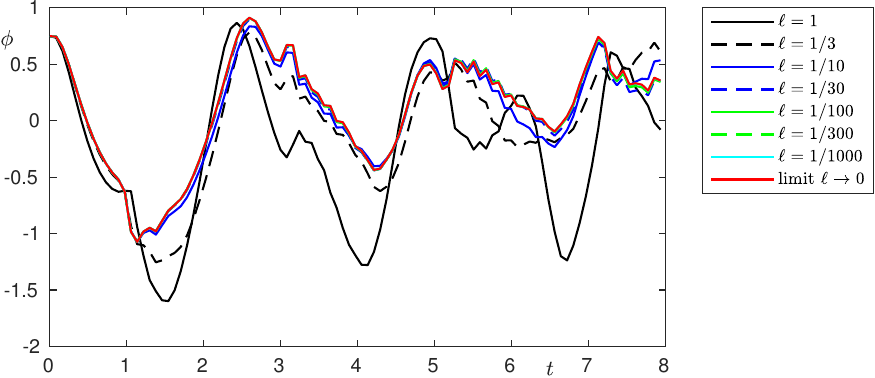}
		\caption{Transversal displacement of the beam, cross-section $x=2$.}
		\label{fig:sl1_first}
	\end{figure}

\begin{figure}[h!]
		\centering
		\includegraphics[width=0.9\textwidth, height=0.27\textheight]{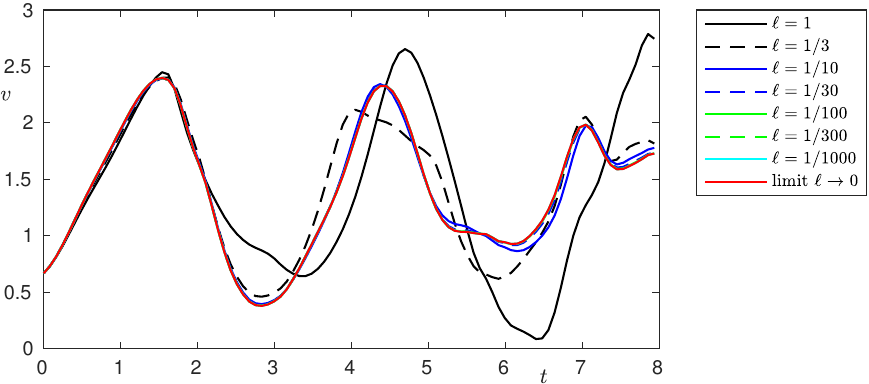}
		\caption{Longitudinal displacement of the beam, cross-section $x=6$.}
		\label{fig:sl1_last}
	\end{figure}

	We perform numerical modelling for the original problem with $l=1,1/3, 1/10, 1/30, 1/100, 1/300, 1/1000$ and the limiting problem ($l=0$) with the following values of constants $\rho_1=\rho_2=1$, $\beta_1=\beta_2=2$, $\sigma_1=4$, $\sigma_2=2$,$\lambda_1=8$, $\lambda_2=4$, $L=10$, $L_0=4$ and the right-hand sides
	\begin{align}
		& p_1(x)=\sin x, && r_1(x)=x, && q_1(x)=\sin x, \label{rhs1}\\
		& p_2(x)=\cos x, && r_2(x)=x+1, && q_2(x)=\cos x. \label{rhs2}
	\end{align}
	In this subsection we consider the nonlinearities with the potentials
	\begin{align*}
		& \F_1(\vph, \psi, \om)=|\vph+\psi|^4-|\vph+\psi|^2+|\vph\psi|^2 + |\om|^3,\\
		& \F_2(u,v,w)=|u+v|^4-|u+v|^2+|uv|^2 + |w|^3.
	\end{align*}
	Consequently, the nonlinearities have the form
	\begin{align*}
		&f_1(\vph, \psi, \om) = 4(\vph+\psi)^3-2(\vph+\psi)+2\vph\psi^2, && f_2(u,v,w)=4(u+v)^3-2(u+v)+2uv^2,\\
		&h_1(\vph, \psi, \om) = 4(\vph+\psi)^3-2(\vph+\psi)+2\vph^2\psi, && h_2(u,v,w)=4(u+v)^3-2(u+v)+2u^2v, \\
		&g_1(\vph, \psi, \om)=3|\om|\om, && g_2(u,v,w)=3|w|w.
	\end{align*}
	For modelling we choose the following (globally Lipschitz) dissipation 
	\begin{equation*}
		\ga(s)=\left\{
		\begin{aligned}
			& \frac{1}{100}s^3, &&|s|\le 10, \\
			& 10s, && |s|>10
		\end{aligned}
		\right.
	\end{equation*}
	
	
	and the following initial data:
	\begin{align*}
		&\vph(x,0)=-\frac{3}{16}x^2+\frac{3}{4}x, && u(x,0)=0,\\
		&\psi(x,0)=-\frac{1}{12}x^2+\frac{7}{12}x, && v(x,0)=-\frac{1}{6}x+\frac{5}{3},\\
		&\om(x,0)=\frac{1}{16}x^2-\frac{1}{4}x, && w(x,0)=-\frac{1}{12}x^2+\frac 76 x-\frac{10}{3},\\
		& \vph_t(x,0)=\frac{x}{4}, && u_t(x,0)=-\frac{1}{6}(x-10), \\
		& \psi_t(x,0)=\frac{x}{4}, && v_t(x,0)=-\frac{1}{6}(x-10), \\
		& \om_t(x,0)=\frac{x}{4}, && w_t(x,0)=-\frac{1}{6}(x-10).
	\end{align*}
	Figures \ref{fig:sl1_first}-\ref{fig:sl1_last} show the behavior of solutions when $l\arr 0$ for the chosen cross-sections of the beam. 
	
	
	\begin{figure}[h!]
		\centering
		\includegraphics[width=0.9\textwidth, height=0.27\textheight]{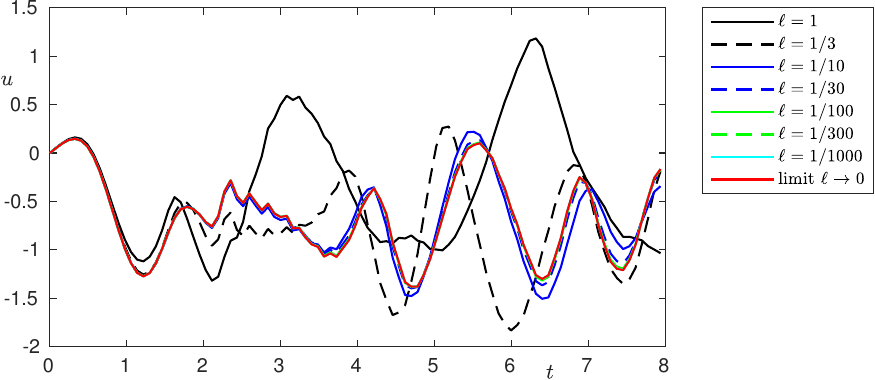}
		\caption{Transversal displacement of the beam, cross-section $x=6$.}
	\end{figure}
	
	\begin{figure}[h!]
		\centering
		\includegraphics[width=0.9\textwidth, height=0.27\textheight]{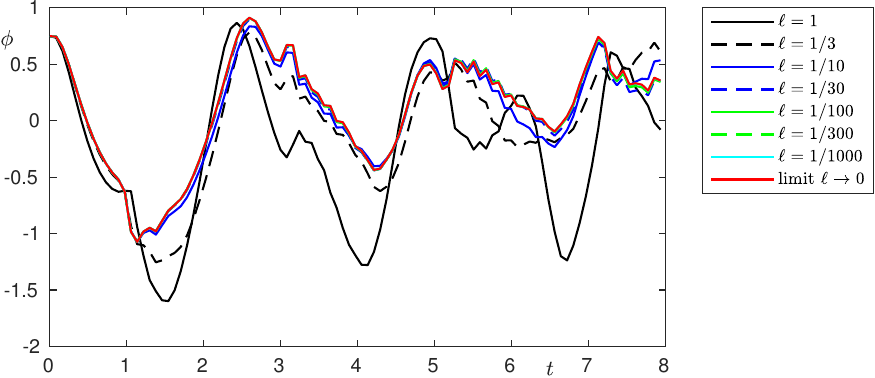}
		\caption{Shear angle variation of the beam, cross-section $x=2$.}
	\end{figure}
	\begin{figure}[h!]
		\centering
		\includegraphics[width=0.9\textwidth, height=0.27\textheight]{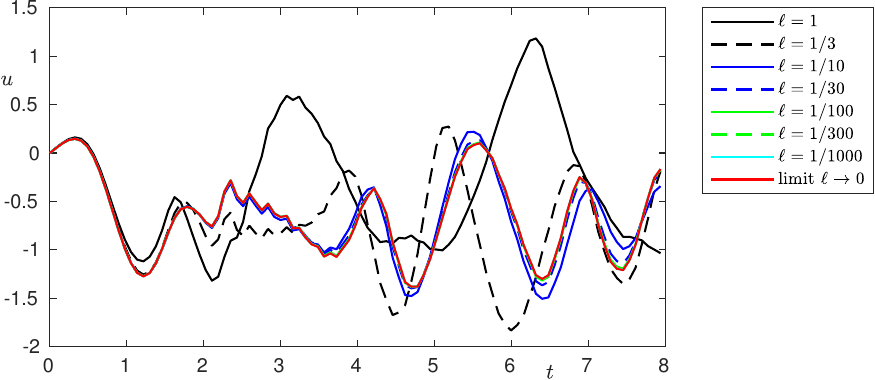}
		\caption{Shear angle variation of the beam, cross-section $x=6$.}
	\end{figure}
	
	\begin{figure}[h!]
		\centering
		\includegraphics[width=0.9\textwidth, height=0.27\textheight]{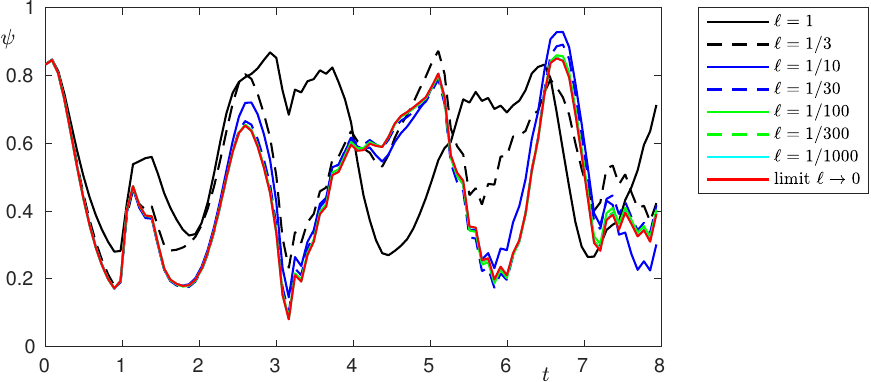}
		\caption{Longitudinal displacement of the beam, cross-section $x=2$.}
	\end{figure}

	\subsection{Singular limit $k_i\arr \infty, \;l\arr 0$}
	
	The singular limit for the straight Timoshenko beam ($l=0$)  as $k_i\arr +\infty$ is the Euler-Bernoulli beam equation \cite[Ch. 4]{Lag1989}. We have a similar result for the Bresse composite beam when $k_i\arr \infty, \;l\arr 0$.
	
	\begin{theorem}
		Let  the assumptions of Theorem \ref{th:WeakWP},  \eqref{Ndecouple} and \eqref{GammaLip1} hold. Moreover,
		\begin{align*}
			\begin{split}
				(\vph_0,u_0)\in\left\{\vph_0\in H^2(0,L_0), \; u_0\in H^2(L_0,L), \; \vph_0(0)=u_0(L)=0,\;\right. \\
				\left. \partial_x\phi_0(0)=\partial_x u_0(L)=0,\; \partial_x\vph_0(L_0,t)=\partial_x u_0(L_0,t) \right\};
			\end{split}\label{ICDsmooth}\tag{I1}\\
			& \psi_0=-\pd_x\vph_{0},\;v_0=-\pd_x u_{0}; \label{ICdep}\tag{I2} \\
			& (\vph_1,u_1) \in\{\vph_1\in H^1(0,L_0), \; u_1\in H^1(L_0,L),\; \vph_1(0)=u_1(L)=0,\; \vph_1(L_0,t)= u_1(L_0,t)\}; \label{ICVSmooth} \tag{I3} \\
			& \omega_0=w_0=0; \label{ICLZero} \tag{I4} \\
			& h_1, h_2\in C^1(\mathbb R^2); \label{NC1}\tag{N6} \\
			\begin{split}
				r_1\in L^\infty(0,T;H^1(0,L_0)), \; r_2\in L^\infty(0,T;H^1(L_0,L)), \\
				r_1(L_0,t)=r_2(L_0,t) \quad \mbox{ for almost all } t>0.
			\end{split}\label{Rsmooth}\tag{R3}
		\end{align*}
		Let  $k_j^{(n)}\arr \infty$, $l^{(n)}\arr 0$ as $n\arr \infty$, and $\Phi^{(n)}$ be weak solutions to \eqref{AEq}-\eqref{AIC} with  fixed $k_j^{(n)}, \; l^{(n)}$ and the same initial data
		\begin{equation*}
			\Phi(x,0)=(\vph_0,\psi_0,\om_0, u_0,v_0,w_0)(x), \quad 	\Phi_t(x,0)=(\vph_1,\psi_1,\om_1, u_1,v_1,w_1).
		\end{equation*}
		Then for every $T>0$
		\begin{align*}
			&\Phi^{(n)}  \stackrel{\ast}{\rightharpoonup} (\vph,\psi,\om, u,v,w) \quad &\mbox{in } L^\infty(0,T;H_d) \; &\mbox{ as } n\arr \infty,\\
			&\Phi^{(n)}_t  \stackrel{\ast}{\rightharpoonup} (\vph_t,\psi_t,\om_t, u_t,v_t,w_t) \quad &\mbox{in } L^\infty(0,T;H_v)\; &\mbox{ as } n\arr \infty,
		\end{align*}
		where
		\begin{itemize}
			\item $(\vph, u)$ is a weak solution to
			\begin{align}
				\begin{split}
					\rho_1\vph_{tt}-\beta_1\vph_{ttxx} +\la_1 \vph_{xxxx} -\ga'(-\vph_{tx}) \vph_{txx} + \pd_x h_1(\vph,-\varphi_x) +{f_1(\vph,-\varphi_x)}= \qquad\\
					p_1(x,{t})+\partial_x r_1(x,{t}), \quad (x,t)\in (0,L_0)\times (0,T), \label{KirchEq1}
				\end{split}\\
				\begin{split}
					\rho_2u_{tt}  -\beta_2u_{ttxx} +\la_2 u_{xxxx}+  \pd_x h_2(u,-u_x)+f_2(u,-u_x)= \qquad \qquad \qquad \qquad\qquad\\
					p_2(x,{t})+\partial_x r_2(x,{t}), \quad (x,t)\in (L_0,L)\times (0,T), \label{KirchEq2}
				\end{split}\\
				&\vph(0,t)=\vph_x(0,t)=0, \, u(L,t)=u_x(L,t)=0,\,\\
				&\vph(L_0,t)=u(L_0,t),   \vph_x(L_0,t)=u_x(L_0,t),
				\la_1 \vph_{xx}(L_0,t)=\la_2 u_{xx}(L_0,t), \, \label{KirchTC1}\\
				\begin{split}
					\la_1 \vph_{xxx}(L_0,t)-\beta_1 \vph_{ttx}(L_0,t)+ h_1(\vph(L_0,t),-\vph_x(L_0,t)) +\ga(-\vph_{tx}(L_0,t))=\qquad\\
					\qquad\la_2 u_{xxx}(L_0,t)-\beta_2 u_{ttx}(L_0,t)+h_2(u(L_0,t),-u_x(L_0,t)),\label{KirchTC2}
				\end{split}			
			\end{align}
			with the initial conditions
			\begin{equation*}
				(\vph, u)(x,0)=(\vph_0, u_0)(x), \quad 	(\vph_t, u_t)(x,0)=(\vph_1, u_1)(x).
			\end{equation*}
			\item $\psi=-\vph_x,  v=-u_x$;
			\item $(\om, w)$ is the solution to
			\begin{align}
				&\rho_1\om_{tt}- \si_1\om_{xx}+{g_1(\om)}=q_1(x,t), \quad (x,t)\in (0,L_0)\times (0,T), \label{KWaveEq}\\
				& \rho_2w_{tt}- \si_2 w_{xx}+g_2(w)=q_2(x,{t}), \quad (x,t)\in (L_0,L)\times (0,T),\\
				&\om(0,t)=0,\,w(L,t)=0,\,\\
				& \si_1 \om_x(L_0,t)=\si_2 w_x(L_0,t), \, \om(L_0,t)=w(L_0,t) \label{KWaveTC}
			\end{align}
			with the initial conditions
			\begin{equation*}
				(\om,w)(x,0)=(0,0), \quad 	(\om_t,w_t)(x,0)=(\om_1, w_1)(x).
			\end{equation*}
		\end{itemize}
	\end{theorem}
	\begin{proof}
		The proof uses the idea from \cite[Ch. 4.3]{Lag1989} and differs from it mainly in transmission conditions. We skip the details of the proof, which coincide with \cite{Lag1989}.\\
		Energy inequality \eqref{EE} implies
		\begin{flalign}
			& \pd_t (\vph^{(n)}, \psi^{(n)}, \om^{(n)}, u^{(n)}, v^{(n)},w^{(n)})  & \mbox{ bounded in } L^\infty(0,T;H_v), \\
			& \psi^{(n)}  & \mbox{ bounded in } L^\infty(0,T;H^1(0,L_0)), \label{ConvF}\\
			& v^{(n)} & \mbox{ bounded in }  L^\infty(0,T;H^1(L_0,L)) \\
			& \om^{(n)}_x - l^{(n)}\vph^{(n)} & \mbox{ bounded in } L^\infty(0,T;L_2(0,L_0)), \\
			& w^{(n)}_x - l^{(n)} u^{(n)} & \mbox{ bounded in } L^\infty(0,T;L_2(L_0,L)), \\
			& k_1^{(n)}(\vph^{(n)}_x + \psi^{(n)} + l^{(n)} \om^{(n)}) & \mbox{ bounded in } L^\infty(0,T;L_2(0,L_0)), \\
			& k_2^{(n)}(u^{(n)}_x + v^{(n)} + l^{(n)} w^{(n)})  & \mbox{ bounded in } L^\infty(0,T;L_2(L_0,L)), \label{ConvL}
		\end{flalign}
		Thus, we can extract subsequences which converge in corresponding spaces weak-$\ast$. Similarly to \cite{Lag1989} we have
		\begin{equation*}
			\vph^{(n)}_x + \psi^{(n)} + l^{(n)} \om^{(n)} \stackrel{\ast}{\rightharpoonup} 0 \quad \mbox{ in }  L^\infty(0,T;L_2(0,L_0)),
		\end{equation*}
		therefore
		\begin{equation*}
			\vph_x =- \psi.
		\end{equation*}
		Analogously,
		\begin{equation*}
			u_x =- v.
		\end{equation*}
		\eqref{ConvF}-\eqref{ConvL} imply
		\begin{align}
			&\om^{(n)} \stackrel{\ast}{\rightharpoonup} \om  &\mbox{ in } L^\infty(0,T;H^1(0,L_0)), &
			&w^{(n)} \stackrel{\ast}{\rightharpoonup} w  &\mbox{ in } L^\infty(0,T;H^1(L_0,L)), \label{Conv1}\\
			&\vph^{(n)} \stackrel{\ast}{\rightharpoonup} \vph  &\mbox{ in } L^\infty(0,T;H^1(0,L_0)), &
			&u^{(n)} \stackrel{\ast}{\rightharpoonup} u  &\mbox{ in } L^\infty(0,T;H^1(L_0,L)). \label{Conv2}
		\end{align}
		Thus, the Aubin's lemma gives that
		\begin{equation}\label{Conv3}
			\Phi^{(n)} \arr \Phi \mbox{ strongly in }  C(0,T; [H^{1-\ep}(0,L_0)]^3\times [H^{1-\ep}(L_0,L)]^3)
		\end{equation}
		for every $\ep>0$ and then
		\begin{equation*}
			\pd_x \vph_0 + \psi_0 + l^{(n)} \om_0 \arr 0 \quad \mbox{ strongly in }  H^{-\ep}(0,L_0),
		\end{equation*}
		This implies that
		\begin{equation*}
			\pd_x \vph_0 =- \psi_0 , \quad \om_0=0.
		\end{equation*}
		Analogously,
		\begin{equation*}
			\pd_x u_0 =- v_0 , \quad w_0=0.
		\end{equation*}
		Let us choose a test function of the form 
		$B=(\be^1,-\be^1_x,0,\be^2,-\be^2_x,0)\in F_T$ such that $\be^1_x(L_0, t)=\be^2_x(L_0, t)$ for almost all $t$.	Due to \eqref{Conv1}-\eqref{Conv3} and Lemma \ref{lem:DissTrans} we can pass to the limit in variational equality \eqref{VEq} as $n\arr\infty$.
		The same way as in \cite[Ch. 4.3]{Lag1989}  we obtain that the limiting functions $\vph, u$ are of higher regularity and satisfy the following variational equality
		\begin{multline}\label{LimVarEq}
			\int_0^T \int_0^{L_0}  \left(\rho_1\vph_t\beta^1_t - \be_1\vph_{tx}\be^1_{tx}\right)dxdt + \int_0^T \int_{L_0}^L \left(\rho_2 u_t\beta^2_t - \be_1 u_{tx}\be^2_{tx}\right)dxdt - \\
			\int_0^{L_0}\left(\rho_1 (\vph_t \beta^1_t)(x,0) - \be_1(\vph_{tx}\be^1_{tx})(x,0)\right)dx + \int_{L_0}^L\left(\rho_2(u_t \beta^2_t)(x,0) - \be_1(u_{tx}\be^2_{tx})(x,0)\right) dx +\\
			\int_0^T \int_0^{L_0}  \la_1\vph_{xx}\be^1_{xx}dxdt  + \int_0^T \int_{L_0}^L\la_2u_{xx} \be^2_{xx}dxdt  - 
			\int_0^T\int_0^{L_0} \ga'(-\vph_{xt})\vph_{txx}\be^1dxdt +\\
			\int_0^T\int_0^{L_0} \left(f_1(\vph,-\vph_x)\be^1 - h_1(\vph,-\vph_x)\be^1_x\right)dxdt + \int_0^T \int_{L_0}^L \left(f_2(u,-u_x)\be^2 - h_2(u,-u_x) \be^2_x\right)dxdt  =\\
			\int_0^T \int_0^{L_0} \left( p_1\be^1 - r_1 \be^1_x\right)dxdt +  \int_0^T \int_{L_0}^L \left(p_2\be^2 - r_2 \be^2_x\right)dxdt.
		\end{multline}
		Provided $\vph, u$ are smooth enough, we can integrate \eqref{LimVarEq} by parts with respect to $x,\; t$ and obtain
		
		\begin{multline}\label{number}
			\int_0^T \int_0^{L_0} (\rho_1-\be_1\pd_{xx})\vph_{tt}\be^1 dxdt + \int_0^T \int_{L_0}^L (\rho_2-\be_2\pd_{xx})u_{tt} \be^2 dxdt + \\
			\int_0^T \left[ \be_1\vph_{ttx}(t,L_0)-\be_2 u_{ttx}(t,L_0)\right] \be^1(t,L_0) dt + \\
			\int_0^T \int_0^{L_0} \la_1\vph_{xxxx}\be^1 dxdt + \int_0^T \int_{L_0}^L \la_2u_{xxxx}\be^2  dxdt +\\
			\int_0^T  \left[\la_1\vph_{xx}-\la_2 u_{xx}\right] (t,L_0) \be^1_x(t,L_0) dt -
			\int_0^T  \left[\la_1\vph_{xxx}-\la_2 u_{xxx}\right] (t,L_0) \be^1(t,L_0) dt  -\\
			\int_0^T \int_0^{L_0} \ga'(-\vph_{xt})\vph_{xxt}\be^1 dxdt - \int_0^T \ga(-\vph_{xt}(L_0,t))\be^1(L_0,t) +\\ 
			\int_0^T\int_0^{L_0} \left(f_1(\vph,-\vph_x) + \pd_x h_1(\vph,-\vph_x)\right)\be^1dxdt   +
			\int_0^T\int_{L_0}^L \left(f_2(u,-u_x) + \pd_x h_2(u,-u_x) \right)\be^2dxdt  + \\
			\int_0^T \left(h_2(u(L_0,t),-u_x(L_0,T)) -h_1(\vph(L_0,t),-\vph_x(L_0,T))\right)\be^1(L_0,t) dt=\\
			\int_0^T \int_0^{L_0} (p_1+\pd_x r_1) \be^1 dxdt + \int_0^T \int_{L_0}^L (p_2+\pd_x r_2) \be^2 dxdt + \int_0^T  \left[ r_2(t,L_0)-r_1(t,L_0)\right] \be^1(t,L_0)dt.
		\end{multline}
		
		Requiring all the terms containing $\be^1(L_0,t)$, $\be^1_x(L_0,t)$ to be zero, we get transmission conditions \eqref{KirchTC1}-\eqref{KirchEq2}. Equations \eqref{KirchEq1}-\eqref{KirchEq2} are recovered from the variational equality \eqref{number}.
		
		Problem \eqref{KWaveEq}-\eqref{KWaveTC} can be obtained in the same way.
	\end{proof}	
	
	We perform numerical modelling for the original problem with the initial parameters
	\begin{equation*}
		l^{(1)}=1, \; k_1^{(1)}=4, \;  k_2^{(1)}=1;\\
	\end{equation*}
	We model the simultaneous convergence $l\arr 0$ and $k_1,k_2\arr \infty$ in the following way: we divide $l$ by the factor $\chi$ and multiply $k_1,k_2$ by the factor $\chi$. Calculations performed for the original problem with 
	\begin{equation*}
		\chi=1,\; \chi=3,\; \chi=10, \; \chi=30, \; \chi=100, \; \chi=300
	\end{equation*}
	and the limiting problem \eqref{KirchEq1}-\eqref{KirchTC2}. 
	Other constants in the original problem are the same as in the previous subsection and we choose the functions in the right-hand side of \eqref{rhs1}-\eqref{rhs2} as follows:
	$$r_1(x)=x+4 ,  \,r_2(x)=2x.$$
	The nonlinear feedbacks are
	\begin{align*}
		& f_1(\vph, \psi,\om)=4\vph^3-2\vph, && f_2(u,v,w)=4u^3-8u, \\
		& h_1(\vph, \psi,\om)=0, && h_2(u,v,w)=0,\\
		& g_1(\vph, \psi,\om)=3|\om|\om, && g_2(u,v,w)=6|w|w.	
	\end{align*}
	We use linear dissipation $\ga(s)=s$ and  choose the following initial displacement and shear angle variation
	$$
	\vph_0(x)=-\frac{13}{640}x^4 +\frac{6}{40}x^2-\frac{23}{40}x^2,
	$$
	$$
	u_0(x)=	\frac{41}{2160}x^4 -\frac{68}{135}x^3 + \frac{823}{180}x^2 -\frac{439}{27}x + \frac{520}{27}.
	$$
	$$
	\psi_0(x)= -\left(-\frac{13}{160}x^3 +\frac{27}{40}x^2 -\frac{23}{20}x\right),
	$$			 
	$$
	v_0(x)=-\left(\frac{41}{540}x^3 -\frac{68}{45}x^2+\frac{823}{90}x -\frac{439}{27}\right).
	$$
	and  set 
	\begin{equation*}
		\om_0(x)=w_0(x)=0.
	\end{equation*}
	We  choose the following initial velocities
	$$
	\vph_1(x)=-\frac{1}{32}x^3+\frac{3}{16}x^2, \quad u_1(x)=\frac{1}{108}x^3 -\frac{7}{36}x^2+\frac{10}{9}x-\frac{25}{27},
	$$
	$$\omega_1(x)=\psi_1(x)=\frac{3}{5}x,$$
	$$w_1(x)=v_1(x)=-\frac{2}{5}x+4.  $$
	
\begin{figure}[h!]
		\centering
		\includegraphics[width=0.9\textwidth, height=0.27\textheight]{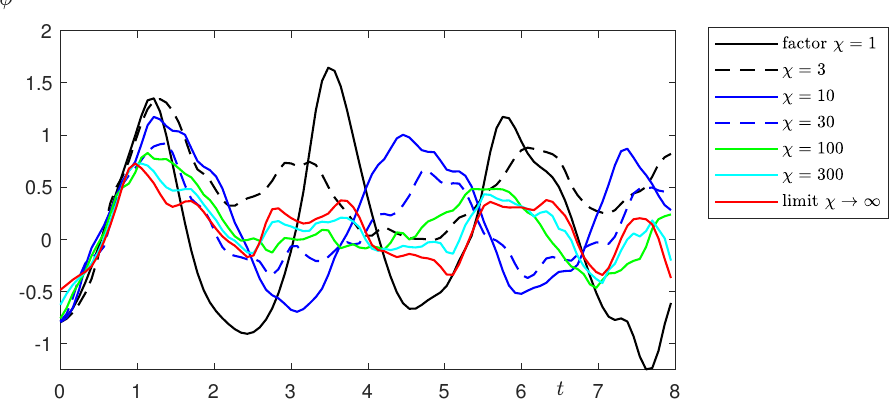}
		\caption{Transversal displacement of the beam, cross-section $x=2$.}
		\label{fig:sl2_first}
	\end{figure}

\begin{figure}[h!]
		\centering
		\includegraphics[width=0.9\textwidth, height=0.27\textheight]{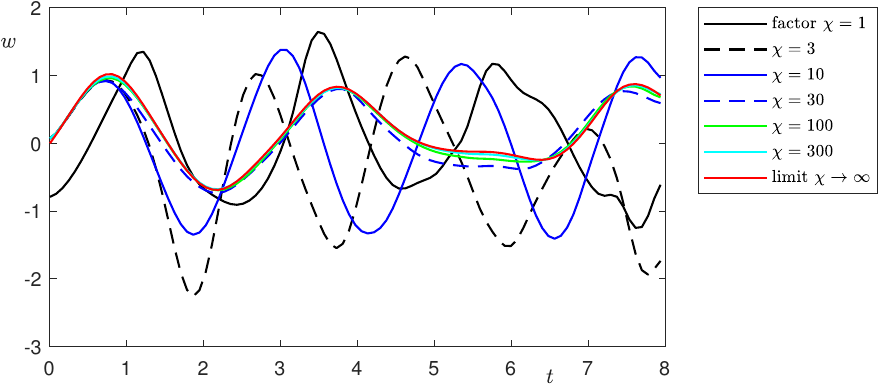}
		\caption{Longitudinal displacement of the beam, cross-section $x=6$.}
		\label{fig:sl2_last}
	\end{figure}

	The double limit case appeared to be more challenging from the  point of view of numerics, then the case $l\arr 0$. The numerical simulations of the coupled system in equations \eqref{Eq1}-\eqref{BC} including the interface conditions in \eqref{TC1}-\eqref{TC4} were done by a semidiscretization of the functions $\phi, \psi,\om, u, v,w$ with respect to the position $x$ and by using an
	explicit scheme for the time integration. That allows to choose the discretized values at grid points near the interface in a separate step so that they obey the transmission conditions. It was necessary to solve a nonlinear system of equations for the six functions at three grid points (at the interface, and left and right
	of the interface) in each time step. Any attempt to use a full implicit numerical scheme led to extremely time-expensive computations due to the large nonlinear system over all discretized values which was to solve in each time step. On the other hand, increasing $k_1,k_2$ increase the stiffness of the system of ordinary differential equations which results from the semidiscretization, and the CFL-conditions requires small time steps ---
	otherwise numerical oscillations occur. Figures~\ref{fig:sl2_first}-\ref{fig:sl2_last} present smoothed numerical solutions, particularly
	necessary for large factors $\chi$, e.\,g. $\chi=300$. When the parameters $k_1,\; k_2$ are large, the material of the beam gets stiff,
	and so does the discretized system of differential  equations. Nevertheless the oscillations are still noticeable in
	the graph. The observation that the factor $\chi$ cannot be arbitrarily enlarged, underlines the importance
	of having the limit problem for $\chi\to\infty$ in \eqref{Eq1}-\eqref{IC}.

	\begin{figure}[h!]
		\centering
		\includegraphics[width=0.9\textwidth, height=0.27\textheight]{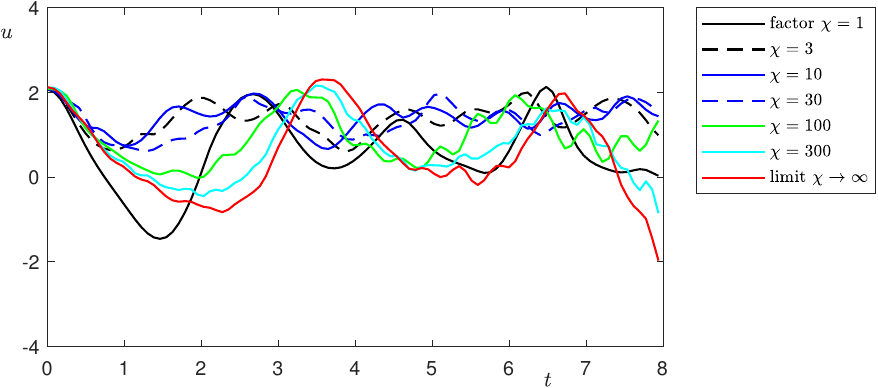}
		\caption{Transversal displacement of the beam, cross-section $x=6$.}
	\end{figure}
	
	\begin{figure}[h!]
		\centering
		\includegraphics[width=0.9\textwidth, height=0.27\textheight]{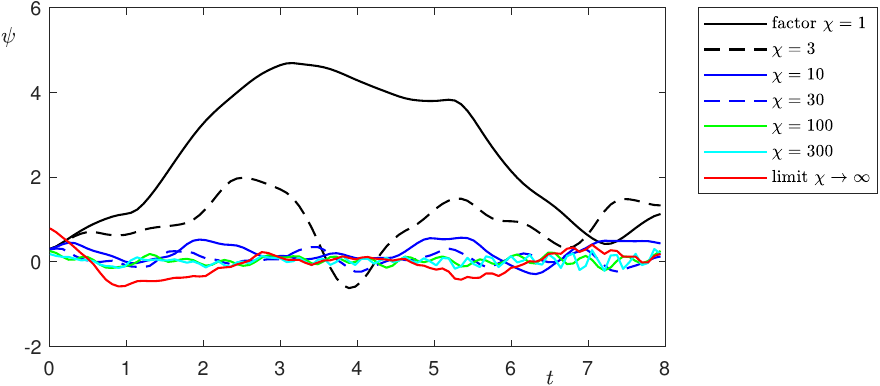}
		\caption{Shear angle variation of the beam, cross-section $x=2$.}
	\end{figure}
	\begin{figure}[h!]
		\centering
		\includegraphics[width=0.9\textwidth, height=0.27\textheight]{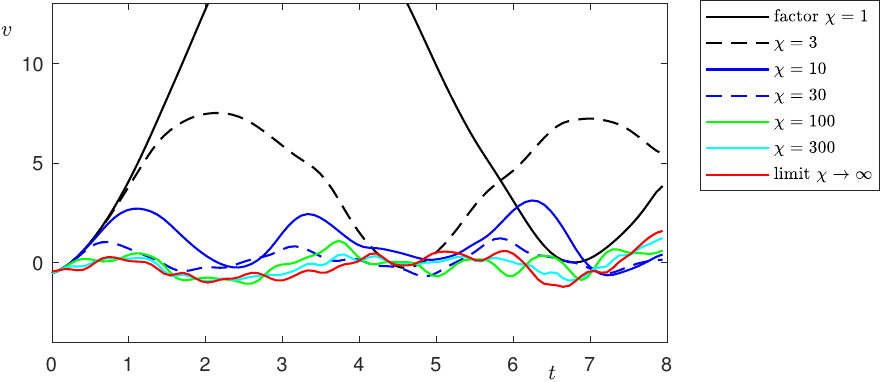}
		\caption{Shear angle variation of the beam, cross-section $x=6$.}
	\end{figure}
	
	\begin{figure}[h!]
		\centering
		\includegraphics[width=0.9\textwidth, height=0.27\textheight]{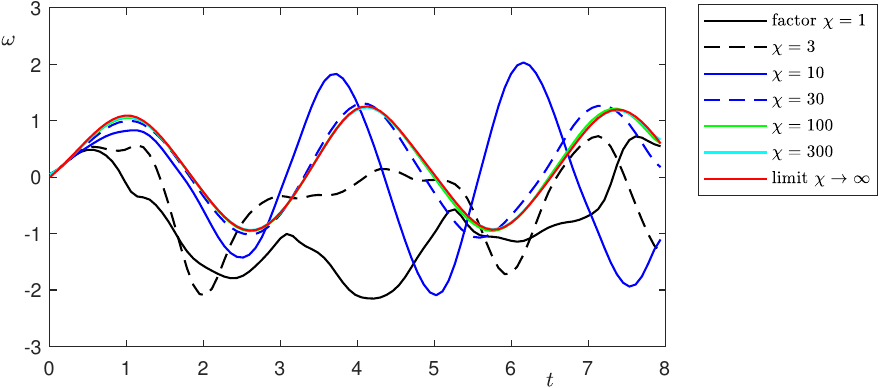}
		\caption{Longitudinal displacement of the beam, cross-section $x=2$.}
	\end{figure}

	\section{Discussion}
	There is a number of papers devoted to long-time behaviour of linear homogeneous Bresse beams (with various boundary conditions and dissipation nature). If damping is present in all three equations, it  appears to be sufficient  for the  exponential stability of the system without additional assumptions on the parameters of the problems (see, e.g., \cite{AlmSan2010}).
	
	The situation is different if we have a dissipation of any kind  in one or two  equations only. First of all, it matters in which equations the dissipation is present. There are  results on the Timoshenko beams \cite{MuRa2002} and the Bresse beams  \cite{Oro2015} that damping in only one of the equations  does not guarantee the exponential stability of the whole system.  It seems that for the Bresse system the presence of the dissipation  in the shear angle equation is necessary for the stability of any kind. To get the exponential stability, one needs  additional assumptions on the coefficients of the problem, usually, the equality of the propagation speeds:
	\begin{equation*}
		k_1=\sigma_1, \quad 	\frac{\rho_1}{k_1}=\frac{\beta_1}{\lambda_1}.
	\end{equation*}
	Otherwise, only polynomial (non-uniform) stability holds (see, e.g., \cite{AlaMuAl2011} for mechanical dissipation and \cite{Oro2015} for thermal dissipation). In  \cite{ChaSoFla2013} analogous results are established in case of nonlinear damping.
	
	If dissipation is present in all three equations of the Bresse system, corresponding problems with nonlinear source forces of local nature possesses  global attractors  under the standard assumptions for nonlinear terms (see, e.g., \cite{MaMo2017}). Otherwise, nonlinear source forces create technical difficulties  and   may cause instability of the system. To the best of our knowledge, there is no literature on such cases.
	
	In the present paper we study a transmission problem for the Bresse system.
	
	Transmission problems for various equation types have already had some history of investigations. One can find a number of papers concerning their well-posedness, long-time behaviour and other aspects (see, e.g., \cite{Pot2012} for a nonlinear thermoelastic/isothermal plate, or \cite{Fast2013} for the Euler-Bernoulli/Timoshenko beam and \cite{Fast2022} for the full von Karman beam). Problems with localized damping are close to transmission problems. In the recent  years a number of such problems for the Bresse beams were studied in, e.g., \cite{MaMo2017,ChaSoFla2013}. To prove the existence of attractors in this case a unique continuation property is an important tool, as well as the frequency method.
	
	The only paper we know on a transmission problem for the Bresse system is \cite{You2022}. The beam in this work consists of a thermoelastic (damped) and elastic (undamped) parts, both purely linear. Despite the presence of dissipation in all three equations for the damped part, the corresponding semigroup is not exponentially stable for any set of parameters, but only polynomially (non-uniformly) stable. In contrast to \cite{You2022}, we consider mechanical damping only in the equation for the shear angle for the damped part. However, we can establish exponential stability for the linear problem and existence of an attractor for the nonlinear one under  restrictions on the coefficients in the damped part only. The assumption on the nonlinearities can be simplified in 1D case (cf. e.g. \cite{Fast2014}).

	\section*{Conflict of Interest Statement}
	The research was conducted in the absence of any commercial or financial relationships that could be construed as a potential conflict of interest.
	
	
	\section*{Acknowledgements}
	The research is supported by the Volkswagen Foundation project "From Modeling and Analysis to Approximation". 
	The first and the third authors were  also successively supported by  the Volkswagen Foundation  project  "Dynamic Phenomena in Elasticity Problems" at Humboldt-Universität zu Berlin, Funding for Refugee Scholars and Scientists from Ukraine.


\end{document}